\documentclass[12pt]{article}
\usepackage{amsfonts}
\usepackage{mathrsfs}
\usepackage{amsmath,amssymb}
\numberwithin{equation}{section}

\newtheorem{theo}{Theorem}[section]
\newtheorem{defi}[theo]{Definition}

\newtheorem{lemm}[theo]{Lemma}

\newtheorem{clai}{Claim}
\newtheorem{conv}{Convention}
\newtheorem{coro}[theo]{Corollary}

\begin{document}\baselineskip 18pt

\centerline{\bf Lie superbialgebra structures on the twisted $N\!=\!1$}
\centerline{{\bf Schr\"{o}dinger-Neveu-Schwarz algebra}}

\vspace*{6pt}

\centerline{Huanxia Fa, Junbo Li}

\centerline{\small School of Mathematics and Statistics, Changshu Institute of Technology, Changshu 215500, China}

\centerline{\small E-mail: sd\_huanxia@163.com, sd\_junbo@163.com}

\vspace*{6pt}

{\small \parskip .005 truein \baselineskip 3pt \lineskip 3pt

\noindent{{\bf Abstract.} Lie superbialgebra
structures on the twisted $N\!=\!1$ Schr\"{o}dinger-Neveu-Schwarz algebra $\frak{tsns}$ are described. The corresponding necessary and sufficient conditions for such superbialgebra to be coboundary triangular are given. Meanwhile, the first cohomology group of $\frak{tsns}$ with coefficients in the tensor product of its adjoint module is completely determined.

\vspace*{6pt}

}

}

\vspace*{18pt}

\leftline{\bf\S1\ \,Preliminaries}
\setcounter{section}{1}\setcounter{theo}{0}\setcounter{equation}{0}

The notion of Lie bialgebras was introduced in 1983 by Drinfeld during the process of investigating quantum groups.
Then there appeared several papers on Lie bialgebras and Lie
superbialgebras (e.g., \cite{M--AM1994,N--JPAA1990,NT--JPAA2000,SS--SCA2006,SS--SCM2015,
YS--AMS2010,YS--CSF2009}). In \cite{M--AM1994,N--JPAA1990,NT--JPAA2000}, the Lie bialgebra structures on Witt and Virasoro algebras were investigated, which are shown to be triangular coboundary and the Lie bialgebra structures on the one-sided Witt algebra were completely classified. In \cite{YS--AMS2010,YS--CSF2009}, the Lie superbialgebra structures on the generalized
super-Virasoro algebra and Ramond $N\!=\!2$ superconformal algebra were investigated. In this paper, we shall study the Lie superbialgebra structures on the twisted $N\!=\!1$ Schr\"{o}dinger-Neveu-Schwarz algebra, which is proved to be coboundary triangular. Both symmetries and super-symmetries act important roles in mathematics and physics. It is known that the Schr\"{o}dinger algebra was realized from the set of dynamic symmetries of the corresponding scalar free Schr\"{o}dinger equation. An $N\!=\!2$ super-symmetric extension of the scalar free Schr\"{o}dinger equation leads to a super-Schr\"{o}dinger model. The Schr\"{o}dinger-Neveu-Schwarz algebras were constructed in  Poisson algebra settings in \cite{HU--NPB2006}, which can be regarded as super-symmetric extensions of the Schr\"{o}dinger algebra.

Firstly, we recall some related definitions. Let
$\mathcal{L}=\mathcal{L}_{\bar0}\oplus\mathcal{L}_{\bar1}$ be a vector space over the complex
number field $\mathbb{C}$. If $x\in\mathcal{L}_{[x]}$, then we say that $x$ is
homogeneous of degree $[x]$ and we write ${\rm deg}x=[x]$. Denote by
$\tau$ the {\it super-twist map} of $\mathcal{L}\otimes\mathcal{L}$:
$\tau(x\otimes y)=(-1)^{[x][y]}y\otimes x$, $\forall\,\,x,y\in\mathcal{L}$. For any $n\in\mathbb{N}$, denote by $\mathcal{L}^{\otimes n}$ the tensor product of
$n$ copies of $\mathcal{L}$ ($\mathcal{L}^{\otimes 2}$ shall be simplified as $\mathcal{L}^{\otimes}$ for convenience) and $\xi$ the {\it super-cyclic map} cyclically
permuting the coordinates of $\mathcal{L}^{\otimes3}$:
$\xi=({\bf\it1}\otimes\tau)\cdot(\tau\otimes{\bf\it1}):\,x_{1}\otimes
x_{2} \otimes x_{3}\mapsto(-1)^{[x_1]([x_2]+[x_3])}x_{2}\otimes
x_{3} \otimes x_{1}$, $\forall\,\,x_i\in\mathcal{L}$, $i=1,2,3$,
where ${\bf\it1}$ is the identity map of $\mathcal{L}$. Then a {\it Lie superalgebra} is a pair $(\mathcal{L},\varphi)$ consisting of a vector space
$\mathcal{L}=\mathcal{L}_{\bar0}\oplus\mathcal{L}_{\bar1}$ and a bilinear map $\varphi:\mathcal{L}\otimes\mathcal{L}\to\mathcal{L}$ satisfying:
\begin{eqnarray*}
&&\varphi(\mathcal{L}_{\bar{i}},\mathcal{L}_{\bar{j}})
\subset\mathcal{L}_{\bar{i}+\bar{j}},\ \
{\rm Ker}({\bf\it1}\otimes{\bf\it1}-\tau)\subset{\rm Ker}\,\varphi,\ \
\varphi\cdot({\bf\it1}\otimes\varphi)\cdot({\bf\it1}
\otimes{\bf\it1}\otimes{\bf\it1}+\xi+\xi^{2})=0.
\end{eqnarray*}
A {\it Lie supercoalgebra} is a pair
$(\mathcal{L},\Delta)$ consisting of a vector space $\mathcal{L}=\mathcal{L}_{\bar0}\oplus
\mathcal{L}_{\bar1}$ and a linear map $\Delta: \mathcal{L} \rightarrow \mathcal{L} \otimes \mathcal{L}$
satisfying:
\begin{eqnarray*}
&&\!\!\!\!\!\!
\Delta(\mathcal{L}_{\bar{i}})\subset
\mbox{$\sum\limits_{\bar{j}\in{\mathbb{Z}_2}}$}
\mathcal{L}_{\bar{j}}\otimes\mathcal{L}_{\bar{i}-\bar{j}},\,\,
{\rm Im}\,\Delta\subset{\rm Im}({\bf\it1}\otimes{\bf\it1}-\tau),\ \,\,
({\bf\it1}\otimes{\bf\it1}\otimes{\bf\it1}+\xi+\xi^{2})
\cdot({\bf\it1}\otimes\Delta)\cdot\Delta=0.
\end{eqnarray*}
A Lie superbialgebra is a triple $(\mathcal{L},\varphi,\Delta)$ satisfying $\Delta\varphi(x\otimes y)=x\ast\Delta y-(-1)^{[x][y]}y\ast\Delta x$, $\forall\,\,x,\,y \in\mathcal{L}$, where $(\mathcal{L},\varphi)$ is a Lie superalgebra and  $(\mathcal{L},\Delta)$ is a Lie super-coalgebra. The symbol ``$\ast$'' means the {\it adjoint diagonal action}:
\begin{eqnarray*}
x\ast(\mbox{$\sum\limits_{i}$}{a_{i}\otimes b_{i}})=
\mbox{$\sum\limits_{i}$}({[x,a_{i}]\otimes b_{i}
+(-1)^{[x][a_i]}a_{i}\otimes[x,b_{i}]}),
\ \ \forall\,\,x,\,a_{i},\,b_{i}\in\mathcal{L}.
\end{eqnarray*}

Denote by $\frak{U}(\mathcal{L})$ the {\it universal enveloping algebra} of
$\mathcal{L}$, ${\bf1}$ the identity element of $\frak{U}(\mathcal{L})$ and $A\backslash B=\{x\,|\,x\in A,x\notin B\}$ for any two sets $A$
and $B$. If $r=\mbox{$\sum\limits_{i}$}{a_{i}\otimes b_{i}}\in\mathcal{L}\otimes\mathcal{L}$, then the following elements are in
$\frak{U}(\mathcal{L})\otimes\frak{U}(\mathcal{L})\otimes\frak{U}(\mathcal{L})$
\begin{eqnarray*}
&&r^{12}=\mbox{$\sum\limits_{i}$}{a_{i}\otimes b_{i}\otimes{\bf1}}=r\otimes{\bf1},\ r^{23} =\mbox{$\sum\limits_{i}$}{{\bf1}\otimes a_{i}\otimes b_{i}}={\bf1}\otimes r,\\
&&r^{13}=\mbox{$\sum\limits_{i}$}{a_{i}\otimes{\bf1}\otimes b_{i}}
=({\bf\it1}\otimes\tau)(r\otimes{\bf1})=(\tau\otimes{\bf\it1})({\bf1}\otimes r),
\end{eqnarray*}
while the following elements are in $\mathcal{L}\otimes\mathcal{L}\otimes\mathcal{L}$
\begin{eqnarray*}
&&[r^{12},r^{23}]=\mbox{$\sum\limits_{i,j}$}a_{i} \otimes[b_{i},a_{j}]\otimes b_{j},\\
&&[r^{12},r^{13}]=\mbox{$\sum\limits_{i,j}$}(-1)^{[a_j][b_i]}
[a_{i},a_{j}]\otimes b_{i}\otimes b_{j},\\
&&[r^{13},r^{23}]=
\mbox{$\sum\limits_{i,j} (-1)^{[a_j][b_i]}$}a_{i}\otimes a_{j} \otimes[b_{i},b_{j}].
\end{eqnarray*}

\begin{defi}\ \ (i) A {\it coboundary superbialgebra} is a quadruple $(\mathcal{L},\varphi,\Delta,r),$
where $(\mathcal{L},\varphi,\Delta)$ is a Lie superbialgebra and
$r\in{\rm Im}({\bf\it1}\otimes{\bf\it1}-\tau)\subset\mathcal{L} \otimes\mathcal{L}$ such that
$\Delta=\Delta_r$ is a {\it coboundary of $r$}, i.e.,
\begin{eqnarray*}
\Delta_r(x)=(-1)^{[r][x]}x\ast r,\ \ \forall\,\,x\in\mathcal{L}.
\end{eqnarray*}
(ii)\ \ A coboundary Lie superbialgebra $(\mathcal{L},\varphi,\Delta,r)$ is
called {\it triangular} if it satisfies the following {\it classical
Yang-Baxter Equation}
\begin{eqnarray}\label{CYBE}
c(r):=[r^{12},r^{13}]+[r^{12},r^{23}]+[r^{13},r^{23}]=0.
\end{eqnarray}
\end{defi}

Let $V=V_{\bar0}\oplus V_{\bar1}$ be an $\mathcal{L}$-module where
$\mathcal{L}=\mathcal{L}_{\bar0}\oplus\mathcal{L}_{\bar1}$. A $\mathbb{Z}_2$-homogenous linear map
$d:\mathcal{L}\to V$ is called a {\it homogenous derivation of degree
$[d]\in\mathbb{Z}_2$}, if $d(\mathcal{L}_i)\subset V_{i+[d]}\ \,(\forall\,\,i\in\mathbb{Z}_2)$,
\begin{eqnarray*}
&&d([x,y])=(-1)^{[d][x]}x\ast d(y)-(-1)^{[y]([d]+[x])}y\ast d(x),
\ \ \forall\,\,x,y\in\mathcal{L}.
\end{eqnarray*}
Denote by ${\rm Der}_{\bar{i}}(\mathcal{L},V)\ \,(\,i=0,1)$ the set of all
homogenous derivations of degree $\bar{i}$. Then the set of all
derivations from $\mathcal{L}$ to $V$
${\rm Der}(\mathcal{L},V)={\rm Der}_{\bar0}(\mathcal{L},V)\oplus{\rm Der}_{\bar1}(\mathcal{L},V)$. Denote
by ${\rm Inn}_{\bar{i}}(\mathcal{L},V)\ \,(\,i=0,1)$ the set of {\it homogenous
inner derivations of degree $\bar{i}$}, consisting of $a_{\rm inn},$
$a\in V_{\bar{i}}$, defined by
\begin{eqnarray*}
&&a_{\rm inn}:\,x\mapsto(-1)^{[a][x]}x\ast a,\ \ \forall\,\,x\in\mathcal{L}.
\end{eqnarray*}
Then the set of inner derivations
${\rm Inn}(\mathcal{L},V)={\rm Inn}_{\bar0}(\mathcal{L},V)\oplus{\rm Inn}_{\bar1}(\mathcal{L},V)$.

Denote by $H^1(\mathcal{L},V)$ the {\it first cohomology group} of $\mathcal{L}$
with coefficients in $V$. Then
\begin{eqnarray*}
H^1(\mathcal{L},V)\cong{\rm Der}(\mathcal{L},V)/{\rm Inn}(\mathcal{L},V).
\end{eqnarray*}
An element $r$ in a superalgebra $\mathcal{L}$ is said to satisfy the
{\it modified Yang-Baxter equation} if
\begin{eqnarray}\label{MYBE}
&&x\ast c(r)=0,\ \ \forall\,\,x\in\mathcal{L}.
\end{eqnarray}

The {\it twisted $N\!=\!1$ Schr\"{o}dinger-Neveu-Schwarz algebra} $\frak{tsns}$ is an
infinite-dimensional Lie superalgebra over the complex field $\mathbb{C}$ with the
basis $\{L_n,G_r,Y_p,M_p\,|\,n\in \mathbb{Z},\,r\in\frac{1}{2}+\mathbb{Z},\,p\in\frac{1}{2}\mathbb{Z}\}$
and the following non-vanishing super brackets:
\begin{eqnarray*}
&&[L_n,L_{m}]=(m-n)L_{m+n},\ \ \ [L_n,G_r\,]=(r-\frac{n}{2})G_{r+n},
\ \ \ \ [G_r,G_s]=2L_{r+s},\\[6pt]
&&[L_n,Y_p]=\left\{\begin{array}{cc}
(p-\frac{n}{2})Y_{p+n}&{\rm if}\ p\in\mathbb{Z},\\[6pt]
pY_{p+n}&{\rm if}\ p\in\frac{1}{2}+\mathbb{Z},
\end{array}\right.
\ \ \,[L_n,M_p]=\left\{\begin{array}{cc}
pM_{p+n}&{\rm if}\ p\in\mathbb{Z},\\[6pt]
(p+\frac{n}{2})M_{p+n}&{\rm if}\ p\in\frac{1}{2}+\mathbb{Z},
\end{array}\right.\\
&&[G_r,Y_p]=
\left\{\begin{array}{cc}
\frac{1}{2}(p-r)Y_{p+r}&{\rm if}\ p\in\mathbb{Z},\\[6pt]
2Y_{p+r}&{\rm if}\ p\in\frac{1}{2}+\mathbb{Z},
\end{array}\right.
\ \,[G_r,M_p]=\left\{\begin{array}{cc}
\frac{p}{2}M_{p+r}&{\rm if}\ p\in\mathbb{Z},\\[6pt]
2M_{p+r}&{\rm if}\ p\in\frac{1}{2}+\mathbb{Z},
\end{array}\right.\\
&&[Y_p,Y_q]=\left\{\begin{array}{ccc}
\frac{1}{2}(q-p)M_{q+p}&{\rm if}\ p,q\in\mathbb{Z},\\[6pt]
\frac{q}{2}M_{q+p}&{\rm if}\ p\in\mathbb{Z},q\in\frac{1}{2}+\mathbb{Z},\\[6pt]
2M_{q+p}&{\rm if}\ p,q\in\frac{1}{2}+\mathbb{Z}.
\end{array}\right.
\end{eqnarray*}
It is easy to see that $\frak{tsns}$ is $\mathbb{Z}_2$-graded with
$\frak{tsns}=\frak{tsns}_{\bar{0}}\oplus\frak{tsns}_{\bar{1}}$, where
\begin{eqnarray*}
\frak{tsns}_{\bar{0}}\!\!\!&=&\!\!\!{\rm{span}}_\mathbb{C}\{L_n,Y_n,M_n\,|\,n\in \mathbb{Z}\},\\
\frak{tsns}_{\bar{1}}\!\!\!&=&\!\!\!{\rm{span}}_\mathbb{C}\{G_r,Y_r,M_r\,|\,r\in\frac{1}{2}+\mathbb{Z}\}.
\end{eqnarray*}
The Cartan subalgebra (exactly the maximal toral subalgebra) of $\frak{tsns}$ is
${\frak{h}}={\mathbb C}L_0\oplus{\mathcal C}$, where ${\mathcal C}={\mathbb C}M_0$
is the center of $\frak{tsns}$. For convenience, we denote by ${\mathcal C}^{\otimes}={\mathbb C}M_0\otimes M_0$. It should be noted that $\frak{tsns}_{\bar{0}}$
is precisely the well-known twisted Schr\"{o}dinger-Virasoro Lie algebra
$\frak{tsv}$ and the subalgebra $\frak{ns}$ spanned by
$\{L_n,G_r\,|\,n\in{\mathbb{Z}},r\in\frac{1}{2}+\mathbb{Z}\}$ is the $N\!=\!1$ Neveu-Schwarz algebra.
For convenience, we denote $\frak{ns}={\rm{span}}_\mathbb{C}\{L_n,G_r\,|\,n\in{\mathbb{Z}},r\in\frac{1}{2}+\mathbb{Z}\}$.
It is easy to see that $\frak{I}={\rm{span}}_\mathbb{C}\{Y_p,M_p\,|\,p\in\frac{1}{2}\mathbb{Z}\}$ is an ideal of
$\frak{tsns}$ and $\frak{tsns}=\frak{ns}\ltimes\frak{I}$.

The following lemma has been obtained in \cite[Theorem 4.2.1]{F--Phd2013}:
\begin{lemm}\label{201303070600}
$H^1(\frak{tsns},\frak{tsns})=\mathfrak{D}$, where the elements of $\mathfrak{D}$ are of the following forms:
\begin{eqnarray*}
&&\mathfrak{d}(L_n)=\alpha nM_{n},\ \ \ \mathfrak{d}(G_s)=s\alpha M_{s},
\ \ \ \mathfrak{d}(Y_p)=\beta Y_{p},\ \ \ \mathfrak{d}(M_p)=2\beta M_{p},
\end{eqnarray*}
\end{lemm}
for any $\alpha,\,\beta \in\mathbb{C}$, $n\in\mathbb{Z}$, $s\in\frac12+\mathbb{Z}$ and $p\in\frac12\mathbb{Z}$.

It is easy to see that $\frak{tsns}$ and $\frak{tsns}^\otimes$ are both $\frac12{\mathbb{Z}}$-graded. Denote ${\rm Der}(\frak{tsns},\frak{tsns}^\otimes)$ \big(resp. ${\rm Inn}(\frak{tsns},\frak{tsns}^\otimes)$\big) the space of derivations (resp. inner derivations) from $\frak{tsns}$ to $\frak{tsns}^\otimes$, and $H^1(\frak{tsns},\frak{tsns}^\otimes)$ the first cohomology group of $\frak{tsns}$ with coefficients in $\frak{tsns}^\otimes$.

The following lemma follows immediately from Lemma \ref{201303070600}.
\begin{lemm}\label{201303070606}
One can find some $\mathfrak{d}^\natural\in{\rm Der}(\frak{tsns},\frak{tsns}^\otimes)$, defined by the following relations:
\begin{eqnarray}
\begin{array}{lll}\label{201303090600}
\mathfrak{d}^\natural(L_n)
\!\!\!&=&\!\!\!\alpha nM_0\otimes M_n+\alpha^\dag nM_n\otimes M_0,\ \ \ \mathfrak{d}^\natural(G_s)=\alpha sM_0\otimes M_s+\alpha^\dag sM_s\otimes M_0,\\[8pt]
\mathfrak{d}^\natural(Y_{p})
\!\!\!&=&\!\!\!\beta M_0\otimes Y_{p}+\beta^\dag Y_p\otimes M_0,
\ \ \ \mathfrak{d}^\natural(M_p)=2(\beta M_0\otimes M_p+\beta^\dag M_p\otimes M_0),
\end{array}
\end{eqnarray}
for any $\alpha,\,\alpha^\dag,\,\beta,\,\beta^\dag\in\mathbb{C}$, $n\in\mathbb{Z}$, $s\in\frac12+\mathbb{Z}$
and $p\in\frac12\mathbb{Z}$.
\end{lemm}

Denote the vector space spanned by $\mathfrak{d}^\natural$ as $\frak{D}^\natural$.
Let $\frak{D}^\natural_0$ be the subspace of $\frak{D}^\natural$ consisting of
elements $\mathfrak{d}^\natural$ such that
$\mathfrak{d}^\natural(\frak{tsns})\subseteq\mathrm{Im}({\bf\it1}\otimes{\bf\it1}-\tau)$.
Namely, $\frak{D}^\natural_0$ is the 2-dimensional subspace of $\frak{D}^\natural$
consisting of elements $\mathfrak{d}^\natural$ with $\alpha=-\alpha^\dag$,
$\beta=-\beta^\dag$.

The main results of this paper can be formulated as the following theorem.
\begin{theo}\label{201303090909}
\begin{itemize}
\item[\rm(i)]
${\rm Der}(\frak{tsns},\frak{tsns}^\otimes)={\rm Inn}(\frak{tsns},\frak{tsns}^\otimes)\oplus\frak{D}^\natural$,
i.e., $H^1(\frak{tsns},\frak{tsns}^\otimes)\cong\frak{D}^\natural$.
\item[\rm(ii)]
Let $(\frak{tsns},[\cdot,\cdot],\Delta)$ be a Lie superbialgebra with
$\Delta=\Delta_r+\mathfrak{d}^\natural$, $r\in\frak{tsns}^\otimes\,({\rm mod\,}{\mathcal C}^\otimes)$ and $\mathfrak{d}^\natural\in\mathfrak{D}^\natural$. Then $r\in\mathrm{Im}({\bf\it1}\otimes{\bf\it1}-\tau)$ and $\mathfrak{d}^\natural\in\mathfrak{D}^\natural_0$. In particular, $(\frak{tsns},[\cdot,\cdot],\mathfrak{d}^\natural)$ is a Lie superbialgebra provided $\mathfrak{d}^\natural\in\mathfrak{D}^\natural_0$.
\item[\rm(iii)]
Let $(\frak{tsns},[\cdot,\cdot],\Delta)$ be a Lie superbialgebra with
$\Delta=\Delta_r+\mathfrak{d}^\natural$, $r\in\frak{tsns}^\otimes\,({\rm mod\,}{\mathcal C}^\otimes)$ and $\mathfrak{d}^\natural\in\mathfrak{D}^\natural_0$. Then it is triangular coboundary if and only if $\alpha=\beta=0$ referred in \eqref{201303090600}.
In other words, $(\frak{tsns},[\cdot,\cdot],\Delta_r+\mathfrak{d}^\natural)$ stands no possibility to be triangular coboundary for any nontrivial $\mathfrak{d}^\natural\in\mathfrak{D}^\natural_0$.
\end{itemize}
\end{theo}

\vspace*{18pt}

\leftline{\bf\S2\ \,Proof of Theorem \ref{201303090909}}
\setcounter{section}{2}\setcounter{theo}{0} \setcounter{equation}{0}

The following result for the non-super case can be found in \cite{NT--JPAA2000}, while its super case can be found in \cite{YS--CSF2009}.
\begin{lemm}\label{201303090606}
Let $\mathcal{L}$ be a Lie superalgebra,
$r\in{\rm Im}({\bf\it1}\otimes{\bf\it1}-\tau)\subset\mathcal{L}\otimes\mathcal{L}$ with $[r]=\bar0$.
Then
\begin{equation}\label{201303090612}
({\bf\it1}\otimes{\bf\it1}\otimes{\bf\it1}+\xi+\xi^{2})\cdot({\bf\it1}\otimes\Delta_r)
\cdot\Delta_r(x)=x\ast c(r),\ \ \forall\,\,x\in\mathcal{L}.
\end{equation}
Thus $(\mathcal{L},[\cdot,\cdot],\Delta_r)$ is a Lie superbialgebra if and only
if $r$ satisfies \eqref{MYBE}.
\end{lemm}

The following lemma can be found in \cite[Lemma 2.2]{LPZ-JA2012}.
\begin{lemm}\label{201303091801}
Suppose that $\frak{g}=\oplus_{n\in\mathbb{Z}}\frak{g}_n$ is a $\mathbb{Z}$-graded Lie algebra with a finite-dimensional center $\frak{C_g}$, and $\frak{g}_0$ is generated by $\{\frak{g}_n,\,n\neq0\}$. Then
\begin{eqnarray*}
H^1(\frak{g},\frak{C_g}\otimes\frak{g}+\frak{g}\otimes\frak{C_g})_0
=\frak{C_g}\otimes H^1(\frak{g},\frak{g})_0+H^1(\frak{g},\frak{g})_0\otimes\frak{C_g}.
\end{eqnarray*}
\end{lemm}
It is not difficult for us to obtain the corresponding result on $\frak{tsns}$.
\begin{lemm}\label{201303091802}
$H^1(\frak{tsns},\mathcal{C}\otimes\frak{tsns}+\frak{tsns}\otimes\mathcal{C})_0
=\mathcal{C}\otimes H^1(\frak{tsns},\frak{tsns})_0+H^1(\frak{tsns},\frak{tsns})_0\otimes\mathcal{C}$.
\end{lemm}

The following lemma has been proved by \cite{Y--FMC2009}.
\begin{lemm}\label{201303091800}
Every Lie superbialgebra structure on the $N\!=\!1$ Neveu-Schwarz algebra $\frak{ns}$ is triangular coboundary and $H^1(\frak{ns},\frak{ns}^\otimes)={\rm Der}(\frak{ns},\frak{ns}^\otimes)/{\rm Inn}(\frak{ns},\frak{ns}^\otimes)=0$.
\end{lemm}

It is known that $\frak{tsns}^{\otimes n}$ can be regarded as a $\frak{tsns}$-module under the adjoint diagonal action of $\frak{tsns}$:
\begin{eqnarray*}
x\ast(v_{1}\otimes v_{2}\otimes\cdots\otimes v_{n})\!\!\!&=&\!\!\!
[x,v_{1}]\otimes v_{2}\otimes\cdots\otimes v_{n}
+(-1)^{[x][v_1]}v_{1}\otimes[x,v_{2}]\otimes\cdots\otimes v_{n}\\
\!\!\!&&\!\!\!+\cdots+(-1)^{[x]([v_1]+\cdots+[v_{n-1}])}v_{1}\otimes v_{2}\otimes\cdots\otimes [x,v_{n}]
\end{eqnarray*}
for all $x,\,v_{i}\in\frak{tsns}$ with $i=1,\,2,\,\cdots,\,n$. The following lemma can be obtained by employing the similar techniques of \cite[Proposition 3.5]{SS--SCA2006} and  \cite[Lemma 2.2]{WSS--CA2007}.
\begin{lemm}\label{201303090800}
If $x\ast r=0$ for any $x\in\frak{tsns}$ and some $r\in\frak{tsns}^{\otimes n}$, then $r\in{\mathcal C}^{\otimes n}$.
\end{lemm}

As a conclusion of Lemma \ref{201303090800}, one immediately obtains the following corollary.
\begin{coro}\label{colo}
An element $r \in {\rm Im}({\bf\it1}\otimes{\bf\it1}-\tau) \subset\frak{tsns}\otimes\frak{tsns}
$ satisfies \eqref{CYBE} if and only if it satisfies \eqref{MYBE}.
\end{coro}

In order to prove Theorem \ref{201303090909}\,(i), we need to make more preparations.

Note that $\frak{tsns}^\otimes=\oplus_{i\in\frac{1}{2}\mathbb{Z}}\frak{tsns}^\otimes_i$ is also $\frac{1}{2}\mathbb{Z}$-graded with $\frak{tsns}^\otimes_i=\sum_{j+k=i} \frak{tsns}_j\otimes\frak{tsns}_k$, where
$i,j,k\in\frac{1}{2}\mathbb{Z}$. We say a derivation $d\in{\rm Der}(\frak{tsns} ,\frak{tsns}^\otimes)$ is
{\it homogeneous of degree $i\in\frac{1}{2}\mathbb{Z}$} if $d(\frak{tsns}^\otimes_j) \subset
\frak{tsns}^\otimes_{i +j}$ for all $j\in\frac{1}{2}\mathbb{Z}$. Set ${\rm Der}(\frak{tsns} , \frak{tsns}^\otimes)_i
=\{d\in {\rm Der}(\frak{tsns} , \frak{tsns}^\otimes) \,|\,{\rm deg\,}d =i\}$ for
$i\in\frac{1}{2}\mathbb{Z}$.

For any $d\in{\rm Der}(\frak{tsns},\frak{tsns}^\otimes)$, $i\in\frac{1}{2}\mathbb{Z}$, $u\in\frak{tsns} _j$ with
$j\in\frac{1}{2}\mathbb{Z}$, we can write $d(u)=\sum_{k\in\frac{1}{2}\mathbb{Z}}v_k\in \frak{tsns}^\otimes$
with $v_k\in \frak{tsns}^\otimes_k$, then we set $d_i(u)=v_{i+j}$. Then
$d_i\in{\rm Der}(\frak{tsns},\frak{tsns}^\otimes)_i$ and
\begin{eqnarray}\label{201303091200}
\mbox{$d=\sum\limits_{i\in\frac{1}{2}\mathbb{Z}} d_i\
\,\mbox{\ where\ }d_i \in {\rm Der}(\frak{tsns}, \frak{tsns}^\otimes)_i,$}
\end{eqnarray}
which holds in the sense that for every $u \in\frak{tsns} $ only finitely
many $d_i(u)\neq 0,$ and $d(u)=\sum_{i \in\frac{1}{2}\mathbb{Z}} d_i(u)$ (we call
such a sum in (\ref{201303091200}) {\it summable}).

Denote $\mathcal{H}=\frak{tsns}\otimes\frak{I}+\frak{I}\otimes\frak{tsns}$. Then $\mathcal{H}$ is a $\frak{tsns}$-submodule of $\frak{tsns}^\otimes$, since $\frak{I}$ is an ideal of $\frak{tsns}$ and denote the quotient $\frak{tsns}$-module $\frak{tsns}^\otimes/\mathcal{H}$ as $\mathcal{Q}$, on which $\frak{I}$ acts trivially and $\mathcal{Q}^\frak{I}=\mathcal{Q}$. The exact sequence $0\rightarrow\mathcal{H}\rightarrow\frak{tsns}^\otimes\rightarrow\frak{tsns}^\otimes/\mathcal{H}\rightarrow0$ induces the following long exact sequence
\begin{eqnarray*}
\longrightarrow H^0(\frak{tsns},\mathcal{Q})\longrightarrow H^1(\frak{tsns},\mathcal{H})\longrightarrow H^1(\frak{tsns},\frak{tsns}^\otimes)\longrightarrow H^1(\frak{tsns},\mathcal{Q})\longrightarrow
\end{eqnarray*}
of $\frac12\mathbb{Z}$-graded vector spaces, where all coefficients of the tensor products are in $\mathbb{C}$. It is easy to see that $H^0(\frak{tsns},\mathcal{Q})=\mathcal{Q}^{\frak{tsns}}=\{x\in\mathcal{Q}\,|\frak{tsns}\ast x=0\,\}=0$. Then
\begin{eqnarray}\label{201303100600}
&&H^1(\frak{tsns},\mathcal{H})\cong H^1(\frak{tsns},\frak{tsns}^\otimes)\ \ \,\mbox{if we can prove}\ \ \, H^1(\frak{tsns},\mathcal{Q})=0.
\end{eqnarray}
Denote $\frak{tsns}_{\mathcal{C}}=\frak{tsns}\otimes{\mathcal{C}}+{\mathcal{C}}\otimes\frak{tsns}$. Then $\frak{tsns}_{\mathcal{C}}$ is an $\frak{tsns}$-submodule of $\mathcal{H}$. The exact sequence $0\rightarrow\frak{tsns}_{\mathcal{C}}\rightarrow\mathcal{H}\rightarrow\mathcal{H}/\frak{tsns}_{\mathcal{C}}\rightarrow0$ induces the following long exact sequence
\begin{eqnarray*}
&&\!\!\!\!\!\!\longrightarrow H^0(\frak{tsns},\mathcal{H}/\frak{tsns}_{\mathcal{C}})\longrightarrow H^1(\frak{tsns},\frak{tsns}_{\mathcal{C}})\longrightarrow H^1(\frak{tsns},\mathcal{H})\longrightarrow H^1(\frak{tsns},\mathcal{H}/\frak{tsns}_{\mathcal{C}})\longrightarrow.
\end{eqnarray*}
It is easy to see that $H^0(\frak{tsns},\mathcal{H}/\frak{tsns}_{\mathcal{C}})=(\mathcal{H}/\frak{tsns}_{\mathcal{C}})^{\frak{tsns}}=\{x\in\mathcal{H}/\frak{tsns}_{\mathcal{C}}\,|\frak{tsns}\cdot x=0\,\}=0$. Then
\begin{eqnarray}\label{201303100606}
&&H^1(\frak{tsns},\frak{tsns}_{\mathcal{C}})\cong H^1(\frak{tsns},\mathcal{H})\ \ \,\mbox{if we can prove}\ \ \,H^1(\frak{tsns},\mathcal{H}/\frak{tsns}_{\mathcal{C}})=0.
\end{eqnarray}

In the following, the notation ``$\,\equiv\,\cdots\,$'' always means
``$\,=\,\cdots\,\big({\rm modulo}\,(\mathcal{C}^\otimes)\big)\,$''.

We shall initiate the proof of Theorem \ref{201303090909} from the first assertion.

\noindent{\it Proof of Theorem \ref{201303090909}\,(i)}\ \ \,It shall follows from a series of claims.
\begin{clai}\label{201302182220}
If $p\in\frac{1}{2}\mathbb{Z}^*$, then $d_p\in{\rm Inn}(\frak{tsns},\frak{tsns}^\otimes)$.
\end{clai}
Denote $u=\frac{1}{p}{\scriptscriptstyle\,}d_{p}(L_0)\in\mathcal{V}_{p}$ where $p\in\frac{1}{2}\mathbb{Z}^*$. For any $x_{q}\in\frak{tsns}_{q}$ with $q\in\frac{1}{2}\mathbb{Z}$, applying $d_{p}$ to
$[L_0,x_{q}]=qx_{q},$ using $d_{p}(x_q)\in\mathcal{V}_{p+q}$ and the
action of $L_0$ on $\mathcal{V}_{p+q}$ is the scalar $p+q$, one has
\begin{eqnarray}\label{201302190808}
(p+q)d_{p}(x_{q})-(-1)^{[d_{p}][x_{q}]}x_{q}\ast d_{p}(L_0)=qd_{p}(x_{q}),
\end{eqnarray}
i.e., $d_{p}(x_{q})=u_{\rm inn}(x_{q})$, which implies $d_{p}$ is
inner. Then this claim follows.

\begin{clai}\label{201302182222}
$d_0(L_0)\equiv d_0(M_0)\equiv0$.
\end{clai}
For any $x\in\frak{tsns}$, taking $p=0$ in \eqref{201302190808}, we obtain $x\ast d_0(L_0)=0$, which together with Lemma \ref{201303090800} gives $d_0(L_0)\equiv0$. For any $x\in\frak{tsns}$, one has $d_0([M_0,x])=0$, which forces $x\ast d_0(M_0)=0$. Then $d_0(M_0)\equiv0$ follows from Lemma \ref{201303090800}.

\begin{clai}\label{201303091202}
$H^1(\frak{tsns},\mathcal{Q})=0$ and $H^1(\frak{tsns},\mathcal{H})\cong H^1(\frak{tsns},\frak{tsns}^\otimes)$.
\end{clai}

The exact sequence $0\rightarrow\frak{I}\rightarrow\frak{tsns}\rightarrow\frak{tsns}/\frak{I}\rightarrow0$ induces an exact sequence of low degree in the Hochschild-Serre spectral sequence
\begin{eqnarray*}
0\longrightarrow H^1(\frak{tsns}/\frak{I},{\mathcal{Q}}^{\frak{I}})\longrightarrow H^1(\frak{tsns},\mathcal{Q})\longrightarrow {H^1(\frak{I},\mathcal{Q})}^{\frak{tsns}/\frak{I}}.
\end{eqnarray*}
According to $\frak{tsns}/\frak{I}\cong\frak{ns}$, $\mathcal{Q}^\frak{I}=\mathcal{Q}$ and $\mathcal{Q}\cong\frak{ns}^\otimes$, Lemma \ref{201303090800} forces $H^1(\frak{tsns}/\frak{I},{\mathcal{Q}}^{\frak{I}})=0$. ${H^1(\frak{I},\mathcal{Q})}^{\frak{tsns}/\frak{I}}$ can be embedded into ${\rm Hom}_{U(\frak{ns})}(\frak{I},\frak{ns}^\otimes)$, which can be easily proved to be zero. Then this claim follows from \eqref{201303100600}.

\begin{clai}\label{201303091203}
$H^1(\frak{tsns},\mathcal{H}/\frak{tsns}_{\mathcal{C}})=0$.
\end{clai}

For any $d_0\in{\rm Der}(\frak{tsns},\mathcal{H}/\frak{tsns}_{\mathcal{C}})$, we can write $d_0(L_{1})$ as follows:
\begin{eqnarray*}
&&\!\!\!\!\!\!\!\!
d_0(L_{1})=\mbox{$\sum\limits_{i\in\mathbb{Z}}$}(a^{L\!Y}_{1,i}L_{i+1}\otimes Y_{-i}
+a^{Y\!L}_{1,i}Y_{i+1}\otimes L_{-i}+a^{L\!M}_{1,i}L_{i+1}\otimes M_{-i}
+a^{M\!L}_{1,i}M_{i+1}\otimes L_{-i})\\
&&\!\!\!\!\!\!\!\!
+\mbox{$\sum\limits_{r\in\frac{1}{2}+\mathbb{Z}}$}(
\!a^{G\!Y}_{1,r}\!G_{r+1}\otimes Y_{-r}+a^{Y\!G}_{1,r}\!Y_{r+1}\otimes G_{-r}
+a^{G\!M}_{1,r}\!G_{r+1}\otimes M_{-r}+a^{M\!G}_{1,r}\!M_{r+1}\otimes G_{-r})\\
&&\!\!\!\!\!\!\!\!
+\mbox{$\sum\limits_{p\in\frac{1}{2}\mathbb{Z}}$}(a^{Y\!Y}_{1,p}Y_{p+1}\otimes Y_{-p}
+a^{Y\!M}_{1,p}Y_{p+1}\otimes M_{-p}+a^{M\!Y}_{1,p}M_{p+1}\otimes Y_{-p}
+a^{M\!M}_{1,p}M_{p+1}\otimes M_{-p}),
\end{eqnarray*}
where the coefficients are all in $\mathbb{C}$ and the sums are all finite.
\begin{conv}\label{201303290600}
The coefficients of $x\otimes M_0$ and $M_0\otimes x$ for all $x\in\frak{tsns}$ should be zero, although we permit them to appear sometimes purely for convenience.
\end{conv}

For any $p\in\frac{1}{2}\mathbb{Z}$ with $x_{p},\,y_{-p}\in\frak{tsns}$, the following identity holds:
\begin{eqnarray*}
L_1\ast({x_{p}\otimes y_{-p}})=
{[L_1,x_{p}]\otimes y_{-p}
+(-1)^{[x_p]}x_{p}\otimes[L_1,y_{-p}]}.
\end{eqnarray*}
Replacing $d_0$ by $d_0 - u_{\rm inn}$, where $u$ is a
combination of some $L_{i}\otimes Y_{-i}$, $Y_{i}\otimes L_{-i}$, $L_{i}\otimes M_{-i}$, $M_{i}\otimes L_{-i}$, $Y_{i}\otimes Y_{-i}$, $Y_{i}\otimes M_{-i}$, $M_{i}\otimes Y_{-i}$, $M_{i}\otimes M_{-i}$, $G_{r}\otimes Y_{-r}$, $Y_{r}\otimes G_{-r}$, $G_{r}\otimes M_{-r}$, $M_{r}\otimes G_{-r}$, $Y_{r}\otimes Y_{-r}$, $Y_{r}\otimes M_{-r}$, $M_{r}\otimes Y_{-r}$ and $M_{r}\otimes M_{-r}$, one can suppose
\begin{eqnarray*}
a^{L\!Y}_{1,i_2}\!\!\!&=&\!\!\!a^{Y\!L}_{1,i_3}=a^{L\!M}_{1,i_4}
=a^{M\!L}_{1,i_5}=a^{G\!Y}_{1,r_2}=\!a^{Y\!G}_{1,r_3}
=a^{G\!M}_{1,r_4}=a^{M\!G}_{1,r_5}\\
\!\!\!&=&\!\!\!a^{Y\!Y}_{1,p_1}=a^{Y\!M}_{1,p_2}
=a^{M\!Y}_{1,p_3}=a^{M\!M}_{1,p_4}=0,
\end{eqnarray*}
for any $i_2\in\mathbb{Z}\backslash\{1\}$, $i_3\in\mathbb{Z}\backslash\{-2\}$, $i_4\in\mathbb{Z}\backslash\{\pm1\}$, $i_5\in\mathbb{Z}\backslash\{-2,\,0\}$, $r_2\in\frac12+\mathbb{Z}\backslash\{\frac12\}$, $r_3\in\frac12+\mathbb{Z}\backslash\{-\frac32\}$, $r_4\in\frac12+\mathbb{Z}\backslash\{\pm\frac12\}$, $r_5\in\frac12+\mathbb{Z}\backslash\{-\frac32,\,-\frac12\}$, $p_1\in\frac12\mathbb{Z}$, $p_2\in\frac12\mathbb{Z}\backslash\{-1,\,-\frac12\}$, $p_3\in\frac12\mathbb{Z}\backslash\{0,\,-\frac12\}$, $p_4\in\frac12\mathbb{Z}\backslash\{-1,\,0,\,-\frac12\}$. Then we can rewrite $d_0(L_{1})$ as follows (just for convenience, we still use the original notations although they have changed):
\begin{eqnarray*}
d_0(L_{1})\!\!\!&=&\!\!\!
a^{L\!Y}_{1,1}L_{2}\otimes Y_{-1}+a^{Y\!L}_{1,-2}Y_{-1}\otimes L_{2}
+a^{L\!M}_{1,-1}L_{0}\otimes M_{1}+a^{L\!M}_{1,1}L_{2}\otimes M_{-1}
+a^{M\!L}_{1,-2}M_{-1}\otimes L_{2}\\
\!\!\!&&\!\!\!+a^{M\!L}_{1,0}M_{1}\otimes L_{0}
+a^{G\!Y}_{1,\frac12}\!G_{\frac32}\otimes Y_{-\frac12}
+a^{Y\!G}_{1,-\frac32}\!Y_{-\frac12}\otimes G_{\frac32}
+a^{G\!M}_{1,-\frac12}\!G_{\frac12}\otimes M_{\frac12}
+a^{G\!M}_{1,\frac12}\!G_{\frac32}\otimes M_{-\frac12}\\
\!\!\!&&\!\!\!+a^{M\!G}_{1,-\frac32}\!M_{-\frac12}\otimes G_{\frac32}
+a^{M\!G}_{1,-\frac12}\!M_{\frac12}\otimes G_{\frac12}
+a^{Y\!M}_{1,-1}Y_{0}\otimes M_{1}
+a^{Y\!M}_{1,-\frac12}Y_{\frac12}\otimes M_{\frac12}
+a^{M\!Y}_{1,0}M_{1}\otimes Y_{0}\\
\!\!\!&&\!\!\!+a^{M\!Y}_{1,-\frac12}M_{\frac12}\otimes Y_{\frac12}
+a^{M\!M}_{1,-\frac12}M_{\frac12}\otimes M_{\frac12}.
\end{eqnarray*}
It should be remarked that although some inner derivations of $L_1$ and $L_{\pm2}$ have been subtracted during the proof of Lemma \ref{201303091800}, we continue to  subtract the above inner derivations of $L_1$ do not impact the proof of Lemma \ref{201303091800} essentially.

For the given $d_0\in{\rm Der}(\frak{tsns},\mathcal{H}/\frak{tsns}_{\mathcal{C}})$, we can write $d_0(M_{\pm\frac12})$ as follows:
\begin{eqnarray*}
&&\!\!\!\!\!\!\!\!d_0(M_{\pm\frac12})=
\mbox{$\sum\limits_{i\in\mathbb{Z}}$}
(\gamma^{L\!Y}_{\pm\frac12,i}L_{i}\otimes Y_{\pm\frac12-i}
+\gamma^{Y\!L}_{\pm\frac12,i}Y_{i\pm\frac12}\otimes L_{-i}
+\gamma^{L\!M}_{\pm\frac12,i}L_{i}\otimes M_{\pm\frac12-i}
+\gamma^{M\!L}_{\pm\frac12,i}M_{i\pm\frac12}\otimes L_{-i})\\
&&\!\!\!\!\!\!\!\!
+\mbox{$\sum\limits_{s\in\frac{1}{2}+\mathbb{Z}}$}
(\gamma^{G\!Y}_{\pm\frac12,s}G_{s}\otimes Y_{\pm\frac12-s}
+\gamma^{Y\!G}_{\pm\frac12,s}Y_{s\pm\frac12}\otimes G_{-s}
+\gamma^{G\!M}_{\pm\frac12,s}G_{s}\otimes M_{\pm\frac12-s}
+\gamma^{M\!G}_{\pm\frac12,s}M_{s\pm\frac12}\otimes G_{-s})\\
&&\!\!\!\!\!\!\!\!
+\mbox{$\sum\limits_{p\in\frac12\mathbb{Z}}$}
(\gamma^{Y\!Y}_{\pm\frac12,p}Y_{p\pm\frac12}\otimes Y_{-p}
+\gamma^{Y\!M}_{\pm\frac12,p}Y_{p\pm\frac12}\otimes M_{-p}
+\gamma^{M\!Y}_{\pm\frac12,p}M_{p\pm\frac12}\otimes Y_{-p}
+\gamma^{M\!M}_{\pm\frac12,p}M_{p\pm\frac12}\otimes M_{-p}),
\end{eqnarray*}
where the coefficients are all in $\mathbb{C}$ and the sums are all finite. The identity $d_0([M_{-\frac12},M_{-\frac12}])=0$ gives $M_{-\frac12}\ast d_0(M_{-\frac12})=0$, which further yields the following identities:
\begin{eqnarray}
\begin{array}{lll}\label{201303290601}
&&(1-i)\gamma^{L\!Y}_{-\frac12,i}=
(1+i)\gamma^{Y\!L}_{-\frac12,i}=
\gamma^{G\!Y}_{-\frac12,s}=\gamma^{Y\!G}_{-\frac12,s}=0,\\[8pt]
&&(i-1)\gamma^{L\!M}_{-\frac12,i}+(i+1)\gamma^{M\!L}_{-\frac12,i}=
\gamma^{G\!M}_{-\frac12,s}+\gamma^{M\!G}_{-\frac12,s}=0,
\end{array}
\end{eqnarray}
for all $i\in\mathbb{Z}$, $s\in\frac12+\mathbb{Z}$. The identities $d_0([G_{\frac12},M_{-\frac12}])=2d_0(M_0)$ and $[L_1,M_{-\frac12}]=0$ yield
\begin{eqnarray}
&&G_{\frac12}\ast d_0(M_{-\frac12})\equiv-M_{-\frac12}\ast d_0(G_{\frac12}),\label{201303290602}\\[8pt]
&&L_1\ast d_0(M_{-\frac12})=M_{-\frac12}\ast d_0(L_1).
\label{201303290603}
\end{eqnarray}
According to the identities given in \eqref{201303290603} and finiteness of the relative sums, we can deduce the following results:
\begin{eqnarray}\label{201303290606}
\gamma^{L\!M}_{-\frac12,i_1}\!\!\!&=&\!\!\!\gamma^{M\!L}_{-\frac12,i_2}=
\gamma^{G\!M}_{-\frac12,s_1}=\gamma^{M\!G}_{-\frac12,s_1}=
\gamma^{Y\!Y}_{-\frac12,p}=\gamma^{Y\!M}_{-\frac12,p}
=\gamma^{M\!Y}_{-\frac12,p}=\gamma^{M\!M}_{-\frac12,p_1}=0,
\end{eqnarray}
for all $i_1\in\mathbb{Z}\backslash\{0,\,1\}$, $i_2\in\mathbb{Z}\backslash\{0,\,-1\}$, $s_1\in\frac12+\mathbb{Z}\backslash\{\pm\frac12\}$, $p\in\frac12\mathbb{Z}$ and $p_1\in\frac12\mathbb{Z}\backslash\{\pm\frac12,\,0,\,1\}$. Combining \eqref{201303290601}, \eqref{201303290602} and \eqref{201303290603}, we also obtain the following identities:
\begin{eqnarray}\label{201303290608}
&&\!\!\!\!\!\!\!\!\!\!\!\!
\gamma^{L\!M}_{-\frac12,0}=-\gamma^{L\!M}_{-\frac12,1}
=\gamma^{M\!L}_{-\frac12,0}=-\gamma^{M\!L}_{-\frac12,-1}=
-\gamma^{G\!M}_{-\frac12,-\frac12}=\gamma^{G\!M}_{-\frac12,\frac12}
=\gamma^{M\!G}_{-\frac12,-\frac12}=-\gamma^{M\!G}_{-\frac12,\frac12}.
\end{eqnarray}
Using the identities related to the coefficients of $d_0(M_{-\frac12})$  given in \eqref{201303290601}, \eqref{201303290602}, \eqref{201303290606} and \eqref{201303290608}, we can simplify $d_0(M_{-\frac12})$ as follows:
\begin{eqnarray*}
d_0(M_{-\frac12})\!\!\!&=&\!\!\!
\gamma^{G\!M}_{-\frac12,-\frac12}(G_{-\frac12}\otimes M_{0}
-M_{0}\otimes\!G_{-\frac12}-L_{0}\otimes M_{-\frac12}
-M_{-\frac12}\otimes L_{0}
+L_{1}\otimes M_{-\frac32}\\
\!\!\!&&\!\!\!
+M_{-\frac32}\otimes L_{1}
-G_{\frac12}\otimes M_{-1}+M_{-1}\otimes G_{\frac12})
+\gamma^{M\!M}_{-\frac12,0}M_{-\frac12}\otimes M_{0}
+\gamma^{M\!M}_{-\frac12,\frac12}M_{0}\otimes M_{-\frac12}\\
\!\!\!&&\!\!\!
+\gamma^{M\!M}_{-\frac12,1}M_{\frac12}\otimes M_{-1}
+\gamma^{M\!M}_{-\frac12,-\frac12}M_{-1}\otimes M_{\frac12}.
\end{eqnarray*}
According to Convention \ref{201303290600}, the following identities hold: $\gamma^{G\!M}_{-\frac12,-\frac12}=\gamma^{M\!M}_{-\frac12,0}
=\gamma^{M\!M}_{-\frac12,\frac12}=0$. Then $d_0(M_{-\frac12})$ can be further simplified as follows:
\begin{eqnarray}\label{201303290609}
d_0(M_{-\frac12})\!\!\!&=&\!\!\!
\gamma^{M\!M}_{-\frac12,1}M_{\frac12}\otimes M_{-1}
+\gamma^{M\!M}_{-\frac12,-\frac12}M_{-1}\otimes M_{\frac12}.
\end{eqnarray}
Furthermore, we can deduce the following identities:
\begin{eqnarray}
\begin{array}{llllll}\label{201303290630}
&&2\gamma^{M\!M}_{-\frac12,-\frac12}=-a^{M\!L}_{1,-2}=-4a^{G\!M}_{1,-\frac12},
\ \ \ 2\gamma^{M\!M}_{-\frac12,1}=-a^{L\!M}_{1,1}=4a^{M\!G}_{1,-\frac12},\\[8pt]
&& a^{Y\!L}_{1,-2}=a^{G\!Y}_{1,\frac12}=a^{Y\!G}_{1,-\frac32}=a^{L\!Y}_{1,1}=
a^{M\!L}_{1,0}+4a^{G\!M}_{1,\frac12}=a^{L\!M}_{1,-1}-4a^{M\!G}_{1,-\frac32}=0.
\end{array}
\end{eqnarray}
Using \eqref{201303290609} and \eqref{201303290602}, we can obtain the following identities:
\begin{eqnarray}
\begin{array}{llllllll}\label{201303290636}
&&\!\!\!\!\!\!\!\!\alpha^{Y\!G}_{\frac12,s}=(i+1)\alpha^{Y\!L}_{\frac12,i}=
\alpha^{G\!Y}_{\frac12,s}=(i-1)\alpha^{L\!Y}_{\frac12,i}=
\alpha^{L\!M}_{\frac12,0}-\gamma^{M\!M}_{-\frac12,-\frac12}=
\gamma^{M\!M}_{-\frac12,1}-\alpha^{M\!L}_{\frac12,0}=0\\[8pt]
&&\!\!\!\!\!\!\!\!
\gamma^{M\!M}_{-\frac12,s-\frac12}+\gamma^{M\!M}_{-\frac12,s}
+\alpha^{G\!M}_{\frac12,s}+\alpha^{M\!G}_{\frac12,s-1}=
i(\gamma^{M\!M}_{-\frac12,i}-\alpha^{M\!L}_{\frac12,i-1})
+(i-1)(\gamma^{M\!M}_{-\frac12,i-\frac12}-\alpha^{L\!M}_{\frac12,i})=0,
\end{array}
\end{eqnarray}
for all $i\in\mathbb{Z}$ and $s\in\frac12+\mathbb{Z}$.

For the given $d_0\in{\rm Der}(\frak{tsns},\mathcal{H}/\frak{tsns}_{\mathcal{C}})$, we can write $d_0(G_{\pm\frac12})$ as follows:
\begin{eqnarray*}
&&\!\!\!\!\!\!\!\!d_0(G_{\pm\frac12})=
\mbox{$\sum\limits_{i\in\mathbb{Z}}$}
(\alpha^{L\!Y}_{\pm\frac12,i}L_{i}\otimes Y_{\pm\frac12-i}
+\alpha^{Y\!L}_{\pm\frac12,i}Y_{i\pm\frac12}\otimes L_{-i}
+\alpha^{L\!M}_{\pm\frac12,i}L_{i}\otimes M_{\pm\frac12-i}
+\alpha^{M\!L}_{\pm\frac12,i}M_{i\pm\frac12}\otimes L_{-i})\\
&&\!\!\!\!\!\!\!\!
+\mbox{$\sum\limits_{s\in\frac{1}{2}+\mathbb{Z}}$}
(\alpha^{G\!Y}_{\pm\frac12,s}G_{s}\otimes Y_{\pm\frac12-s}
+\alpha^{Y\!G}_{\pm\frac12,s}Y_{s\pm\frac12}\otimes G_{-s}
+\alpha^{G\!M}_{\pm\frac12,s}G_{s}\otimes M_{\pm\frac12-s}
+\alpha^{M\!G}_{\pm\frac12,s}M_{s\pm\frac12}\otimes G_{-s})\\
&&\!\!\!\!\!\!\!\!
+\mbox{$\sum\limits_{p\in\frac12\mathbb{Z}}$}
(\alpha^{Y\!Y}_{\pm\frac12,p}Y_{p\pm\frac12}\otimes Y_{-p}
+\alpha^{Y\!M}_{\pm\frac12,p}Y_{p\pm\frac12}\otimes M_{-p}
+\alpha^{M\!Y}_{\pm\frac12,p}M_{p\pm\frac12}\otimes Y_{-p}
+\alpha^{M\!M}_{\pm\frac12,p}M_{p\pm\frac12}\otimes M_{-p}),
\end{eqnarray*}
where the coefficients are all in $\mathbb{C}$ and the sums are all finite. The identity $[G_{\frac12},G_{\frac12}]=2L_1$ gives. Using \eqref{201303290636} and comparing the relative coefficients, we can obtain the following identities:
\begin{eqnarray}
&&G_{\frac12}\ast d_0(G_{\frac12})=d_0(L_1).\label{201303300600}
\end{eqnarray}
Recalling \eqref{201303290602} and \eqref{201303290609}, we know that $\alpha^{G\!M}_{\frac12,\frac12}=-\alpha^{M\!G}_{\frac12,-\frac12}$. Furthermore, according to Convention \ref{201303290600}, we also know that
\begin{eqnarray}\label{201303300603}
&&\alpha^{G\!M}_{\frac12,\frac12}=\alpha^{M\!G}_{\frac12,-\frac12}=0.
\end{eqnarray}

Using \eqref{201303290636}, \eqref{201303300600} and \eqref{201303300603}, we can rewrite $d_0(G_{\frac12})$ as follows:
\begin{eqnarray}\label{201303300606}
\!\!\!\!\!\!
d_0(G_{\frac12})\!\!\!&=&\!\!\!
\alpha^{L\!M}_{\frac12,0}L_{0}\otimes M_{\frac12}
+\alpha^{M\!L}_{\frac12,0}M_{\frac12}\otimes L_{0}
+\alpha^{G\!M}_{\frac12,\frac32}G_{\frac32}\otimes M_{-1}\nonumber\\
\!\!\!&&\!\!\!
+\alpha^{M\!G}_{\frac12,-\frac32}M_{-1}\otimes G_{\frac32}
+\alpha^{M\!Y}_{\frac12,0}M_{\frac12}\otimes Y_{0}
+\alpha^{Y\!M}_{\frac12,-\frac12}Y_{0}\otimes M_{\frac12}.
\end{eqnarray}

Combining \eqref{201303290636}, \eqref{201303300600} and \eqref{201303300603}, we can also deduce the following identities:
\begin{eqnarray}
\begin{array}{llllllll}\label{201303300602}
&&2\alpha^{L\!M}_{\frac12,0}=-4a^{G\!M}_{1,-\frac12}=a^{L\!M}_{1,-1},\ \ \
2\alpha^{G\!M}_{\frac12,\frac32}=4a^{G\!M}_{1,\frac12}=a^{L\!M}_{1,1},\\[8pt]
&&2\alpha^{Y\!M}_{\frac12,-\frac12}=a^{Y\!M}_{1,-1}
=-8a^{Y\!M}_{1,-\frac12},\ \ \ 2\alpha^{M\!Y}_{\frac12,0}=a^{M\!Y}_{1,0}
=8a^{M\!Y}_{1,-\frac12},\\[8pt]
&&2\alpha^{M\!G}_{\frac12,-\frac32}=a^{M\!L}_{1,-2}
=-4a^{M\!G}_{1,-\frac32},\ \ \ 2\alpha^{M\!L}_{\frac12,0}=a^{M\!L}_{1,0}=4a^{M\!G}_{1,-\frac12}.
\end{array}
\end{eqnarray}

The identity $[G_{\frac12},G_{-\frac12}]=2L_0$, together with Claim \ref{201302182222}, gives $G_{\frac12}\ast d_0(G_{-\frac12})\equiv -G_{-\frac12}\ast d_0(G_{\frac12})$. Using \eqref{201303300606}, and comparing the relative coefficients, we can deduce the following identities:
\begin{eqnarray}
\begin{array}{llllllll}\label{201303300800}
&&\alpha^{G\!Y}_{-\frac12,s}=\alpha^{L\!Y}_{-\frac12,i}
=\alpha^{Y\!L}_{-\frac12,i}=\alpha^{Y\!G}_{-\frac12,s}=
\alpha^{Y\!Y}_{-\frac12,p}=0,\\[8pt]
&&\alpha^{Y\!M}_{-\frac12,p_1}=\alpha^{M\!Y}_{-\frac12,p_2}
=\alpha^{M\!M}_{-\frac12,p_3}=\alpha^{L\!M}_{-\frac12,i_1}
=\alpha^{G\!M}_{-\frac12,s_1}=\alpha^{M\!L}_{-\frac12,i_2}
=\alpha^{M\!G}_{-\frac12,s_1}=0,
\end{array}
\end{eqnarray}
for all $i\in\mathbb{Z}$, $i_1\in\mathbb{Z}\backslash\{0,\,1\}$, $i_2\in\mathbb{Z}\backslash\{0,\,-1\}$, $s\in\frac12+\mathbb{Z}$, $s_1\in\frac12+\mathbb{Z}\backslash\{\pm\frac12\}$, $p\in\frac12\mathbb{Z}$, $p_1\in\frac12\mathbb{Z}\backslash\{0\}$, $p_2\in\frac12\mathbb{Z}\backslash\{\frac12\}$ and $p_3\in\frac12\mathbb{Z}\backslash\{0,\,\frac12\}$.

According to Convention \ref{201303290600}, the following identities hold:  \begin{eqnarray}\label{201303300802}
&&\alpha^{M\!M}_{-\frac12,0}=\alpha^{M\!M}_{-\frac12,\frac12}
=\alpha^{Y\!M}_{-\frac12,0}=\alpha^{M\!Y}_{-\frac12,\frac12}
=\alpha^{G\!M}_{-\frac12,-\frac12}=\alpha^{M\!G}_{-\frac12,\frac12}=0.
\end{eqnarray}
During the process of comparing the relative coefficients of $G_{\frac12}\ast d_0(G_{-\frac12})\equiv -G_{-\frac12}\ast d_0(G_{\frac12})$, we also can obtain the following identities (together with \eqref{201303300800} and \eqref{201303300802}):
\begin{eqnarray}
\begin{array}{llllllll}\label{201303300801}
\alpha^{M\!G}_{\frac12,-\frac32}
\!\!\!&=&\!\!\!\alpha^{Y\!M}_{\frac12,-\frac12}
=\alpha^{M\!L}_{\frac12,0}=\alpha^{L\!M}_{\frac12,1}
=\alpha^{M\!Y}_{\frac12,0}=\alpha^{G\!M}_{\frac12,\frac32}
=\alpha^{L\!M}_{\frac12,0}\\[8pt]
\!\!\!&=&\!\!\!\alpha^{Y\!M}_{-\frac12,0}=\alpha^{M\!L}_{-\frac12,0}
=\alpha^{M\!L}_{-\frac12,1}=\alpha^{M\!G}_{-\frac12,-\frac32}
=\alpha^{M\!Y}_{-\frac12,\frac12}=\alpha^{L\!M}_{-\frac12,0}
=\alpha^{L\!M}_{-\frac12,1}+\alpha^{G\!M}_{-\frac12,\frac12}=0.
\end{array}
\end{eqnarray}
Then using \eqref{201303300602} and \eqref{201303300801}, one has the following identities:
\begin{eqnarray}
a^{M\!M}_{1,-\frac12}\!\!\!&=&\!\!\!a^{G\!M}_{1,-\frac12}=a^{L\!M}_{1,-1}
=a^{G\!M}_{1,\frac12}=a^{L\!M}_{1,1}=a^{Y\!M}_{1,-1}
=a^{Y\!M}_{1,-\frac12}\nonumber\\
\!\!\!&=&\!\!\!a^{M\!Y}_{1,0}=a^{M\!Y}_{1,-\frac12}=a^{M\!L}_{1,-2}
=a^{M\!G}_{1,-\frac32}= a^{M\!L}_{1,0}
=a^{M\!G}_{1,-\frac12}=0.\label{201303300802}
\end{eqnarray}
Then the following result follows from \eqref{201303290630} and \eqref{201303300802}:
\begin{eqnarray}\label{201303300803}
d_0(L_1)\!\!\!&=&\!\!\!0.
\end{eqnarray}
Combining \eqref{201303290602}, \eqref{201303290603} and \eqref{201303300803}, we can deduce
\begin{eqnarray*}
&&\gamma^{M\!M}_{-\frac12,1}=2a^{M\!G}_{1,-\frac12}=0
=-2a^{G\!M}_{1,-\frac12}=\gamma^{M\!M}_{-\frac12,-\frac12},
\end{eqnarray*}
according to which, we can simplified $d_0(M_{-\frac12})$ referred in \eqref{201303290609} as follows:
\begin{eqnarray}\label{201303300816}
d_0(M_{-\frac12})\!\!\!&=&\!\!\!0.
\end{eqnarray}
According to \eqref{201303300801}, $d_0(G_{\frac12})$ referred in \eqref{201303300606} can be simplified as follows:
\begin{eqnarray}\label{201303300806}
d_0(G_{\frac12})\!\!\!&=&\!\!\!0.
\end{eqnarray}

The identity $[L_1,G_{-\frac12}]=-G_{\frac12}$, together with \eqref{201303300803} and \eqref{201303300806}, gives $L_1\ast d_0(G_{-\frac12})=0$. Comparing the relative coefficients, we can obtain the following identities:
\begin{eqnarray*}
&&\alpha^{L\!M}_{-\frac12,0}+\alpha^{L\!M}_{-\frac12,1}
=\alpha^{M\!L}_{-\frac12,-1}+\alpha^{M\!L}_{-\frac12,0}
=\alpha^{G\!M}_{-\frac12,-\frac12}+\alpha^{G\!M}_{-\frac12,\frac12}
=\alpha^{M\!G}_{-\frac12,-\frac12}+\alpha^{M\!G}_{-\frac12,\frac12}=0,
\end{eqnarray*}
which together with \eqref{201303300802} and \eqref{201303300801}, force
\begin{eqnarray}\label{201303300809}
&&\alpha^{L\!M}_{-\frac12,1}=\alpha^{M\!L}_{-\frac12,-1}
=\alpha^{G\!M}_{-\frac12,\frac12}=\alpha^{M\!G}_{-\frac12,-\frac12}=0.
\end{eqnarray}
Then the following result follows from \eqref{201303300800}, \eqref{201303300802}, \eqref{201303300801} and \eqref{201303300809}:
\begin{eqnarray}\label{201303300810}
d_0(G_{-\frac12})\!\!\!&=&\!\!\!0.
\end{eqnarray}

The identity $[G_{-\frac12},M_{\frac12}]=2M_0$, together with Claim \ref{201302182222} and \eqref{201303300810}, gives $G_{-\frac12}\ast d_0(M_{\frac12})\equiv0$. Comparing the relative coefficients, we can deduce the following identities:
\begin{eqnarray*}
&&\gamma^{Y\!Y}_{\frac12,p}=\gamma^{Y\!M}_{\frac12,p}
=\gamma^{M\!Y}_{\frac12,p}=\gamma^{M\!M}_{\frac12,p}=0,
\ \ \ \forall\,\,p\in\frac12\mathbb{Z}.
\end{eqnarray*}

The identity $[M_{\frac12},M_{\frac12}]=0$ gives $M_{\frac12}\ast d_0(M_{\frac12})=0$. Comparing the relative coefficients, we can obtain the following identities:
\begin{eqnarray}
\gamma^{G\!Y}_{\frac12,s}\!\!\!&=&\!\!\!
\gamma^{Y\!G}_{\frac12,s}=(i+1)\gamma^{L\!Y}_{\frac12,i}=
(i-1)\gamma^{Y\!L}_{\frac12,i}\nonumber\\
\!\!\!&=&\!\!\!(i+1)\gamma^{L\!M}_{\frac12,i}+(i-1)\gamma^{M\!L}_{\frac12,i}=
\gamma^{G\!M}_{\frac12,s}+\gamma^{M\!G}_{\frac12,s}=0,
\ \ \ \forall\,\,i\in\mathbb{Z},\,s\in\frac12+\mathbb{Z}.\label{201303300900}
\end{eqnarray}

The identity $[M_{\frac12},M_{-\frac12}]=0$, together with Claim \ref{201302182222} and \eqref{201303300816}, gives $M_{-\frac12}\ast d_0(M_{\frac12})\equiv0$. Comparing the relative coefficients, we can obtain the following identities:
\begin{eqnarray}
(i-1)\gamma^{L\!Y}_{\frac12,i}\!\!\!&=&\!\!\!i\gamma^{Y\!L}_{\frac12,i-1}=
(i-1)\gamma^{L\!M}_{\frac12,i}+i\gamma^{M\!L}_{\frac12,i-1}\nonumber\\
\!\!\!&=&\!\!\!\gamma^{G\!Y}_{\frac12,s}=\gamma^{Y\!G}_{\frac12,s-1}=
\gamma^{G\!M}_{\frac12,s}+\gamma^{M\!G}_{\frac12,s-1}=0,
\ \ \ \forall\,\,i\in\mathbb{Z},\,s\in\frac12+\mathbb{Z}.\label{201303300901}
\end{eqnarray}
Combining \eqref{201303300900} and \eqref{201303300901}, we can deduce the following identities:
\begin{eqnarray}\label{201303300902}
\gamma^{G\!M}_{\frac12,s}\!\!\!&=&\!\!\!\gamma^{M\!G}_{\frac12,s}
=\gamma^{L\!Y}_{\frac12,i}=\gamma^{Y\!L}_{\frac12,i}
=\gamma^{L\!M}_{\frac12,i}=\gamma^{M\!L}_{\frac12,i}
=\gamma^{G\!Y}_{\frac12,s}=\gamma^{Y\!G}_{\frac12,s}=0,
\end{eqnarray}
for all $i\in\mathbb{Z}$, $s\in\frac12+\mathbb{Z}$. Then the following result follows from \eqref{201303300900} and \eqref{201303300902}:
\begin{eqnarray}\label{201303300906}
d_0(M_{\frac12})\!\!\!&=&\!\!\!0.
\end{eqnarray}

For the given $d_0\in{\rm Der}(\frak{tsns},\mathcal{H}/\frak{tsns}_{\mathcal{C}})$, we can write $d_0(G_{\pm\frac32})$ as follows:
\begin{eqnarray*}
&&\!\!\!\!\!\!\!\!d_0(G_{\pm\frac32})=
\mbox{$\sum\limits_{i\in\mathbb{Z}}$}
(\alpha^{L\!Y}_{\pm\frac32,i}L_{i}\otimes Y_{\pm\frac32-i}
+\alpha^{Y\!L}_{\pm\frac32,i}Y_{i\pm\frac32}\otimes L_{-i}
+\alpha^{L\!M}_{\pm\frac32,i}L_{i}\otimes M_{\pm\frac32-i}
+\alpha^{M\!L}_{\pm\frac32,i}M_{i\pm\frac32}\otimes L_{-i})\\
&&\!\!\!\!\!\!\!\!
+\mbox{$\sum\limits_{s\in\frac{1}{2}+\mathbb{Z}}$}
(\!\alpha^{G\!Y}_{\pm\frac32,s}G_{s}\otimes Y_{\pm\frac32-s}+\alpha^{Y\!G}_{\pm\frac32,s}\!Y_{s\pm\frac32}\otimes G_{-s}
+\alpha^{G\!M}_{\pm\frac32,s}G_{s}\otimes M_{\pm\frac32-s}
+\alpha^{M\!G}_{\pm\frac32,s}M_{s\pm\frac32}\otimes G_{-s})\\
&&\!\!\!\!\!\!\!\!
+\mbox{$\sum\limits_{p\in\frac12\mathbb{Z}}$}
(\alpha^{Y\!Y}_{\pm\frac32,p}Y_{p\pm\frac32}\otimes Y_{-p}
+\alpha^{Y\!M}_{\pm\frac32,p}Y_{p\pm\frac32}\otimes M_{-p}
+\alpha^{M\!Y}_{\pm\frac32,p}M_{p\pm\frac32}\otimes Y_{-p}
+\alpha^{M\!M}_{\pm\frac32,p}M_{p\pm\frac32}\otimes M_{-p}),
\end{eqnarray*}
where the coefficients are all in $\mathbb{C}$ and the sums are all finite. The identity $[L_1,G_{-\frac32}]=-2G_{-\frac12}$, together with \eqref{201303300803} and \eqref{201303300810}, gives $L_1\ast d_0(G_{-\frac32})=0$. Comparing the relative coefficients, we can deduce the following identities:
\begin{eqnarray}
\!\!\!\!\!\!\!\!
\alpha^{L\!Y}_{-\frac32,i}\!\!\!&=&\!\!\!\alpha^{L\!M}_{-\frac32,i_1}
=\alpha^{Y\!L}_{-\frac32,i}=\alpha^{M\!L}_{-\frac32,i_1}
=\alpha^{G\!Y}_{-\frac32,s}=\alpha^{G\!M}_{-\frac32,s_1}
=\alpha^{Y\!G}_{-\frac32,s}=\alpha^{M\!G}_{-\frac32,s_2}\nonumber\\
\!\!\!&=&\!\!\!\alpha^{Y\!Y}_{-\frac32,p}
=\alpha^{Y\!M}_{-\frac32,p}=\alpha^{M\!Y}_{-\frac32,p}
=\alpha^{M\!M}_{-\frac32,p_1}=0,\label{201303310600}
\end{eqnarray}
for all $i\in\mathbb{Z}$, $i_1\in\mathbb{Z}\backslash\{-1,\,0,\,1\}$, $s\in\frac12+\mathbb{Z}$, $s_1\in\frac12+\mathbb{Z}\backslash\{-\frac32,\,-\frac12,\,\frac12\}$, $s_2\in\frac12+\mathbb{Z}\backslash\{-\frac12,\,\frac12,\,\frac32\}$, $p\in\frac12\mathbb{Z}$, $p_1\in\frac12\mathbb{Z}\backslash\{0,\,1,\,\frac12,\,\frac32\}$. During the process of comparing the relative coefficients of $L_1\ast d_0(G_{-\frac32})=0$, we also can obtain the following identities:
\begin{eqnarray}
\begin{array}{llllllll}\label{201303310601}
&&\alpha^{M\!M}_{-\frac32,\frac12}+\alpha^{M\!M}_{-\frac32,\frac32}
=\alpha^{M\!M}_{-\frac32,0}+\alpha^{M\!M}_{-\frac32,1}=0,\\[8pt]
&&\alpha^{L\!M}_{-\frac32,0}=-2\alpha^{L\!M}_{-\frac32,1}
=-2\alpha^{L\!M}_{-\frac32,-1},\ \ \ \ \ \
\alpha^{M\!L}_{-\frac32,0}=-2\alpha^{M\!L}_{-\frac32,1}
=-2\alpha^{M\!L}_{-\frac32,-1},\\[8pt]
&&\alpha^{G\!M}_{-\frac32,-\frac12}=-2\alpha^{G\!M}_{-\frac32,-\frac32}
=-2\alpha^{G\!M}_{-\frac32,\frac12},\ \ \ \ \ \
\alpha^{M\!G}_{-\frac32,\frac12}=-2\alpha^{M\!G}_{-\frac32,-\frac12}
=-2\alpha^{M\!G}_{-\frac32,\frac32}.
\end{array}
\end{eqnarray}
According to Convention \ref{201303290600}, the following identities hold: $\alpha^{G\!M}_{-\frac32,-\frac32}=\alpha^{M\!G}_{-\frac32,\frac32}
=\alpha^{M\!M}_{-\frac32,0}=\alpha^{M\!M}_{-\frac32,\frac32}=0$. Combining \eqref{201303310600} and \eqref{201303310601}, we can rewrite $d_0(G_{\frac32})$ as follows:
\begin{eqnarray}
d_0(G_{-\frac32})\!\!\!&=&\!\!\!
\alpha^{L\!M}_{-\frac32,-1}(L_{-1}\otimes M_{-\frac12}
-2L_{0}\otimes M_{-\frac32}+L_{1}\otimes M_{-\frac52})\nonumber\\
\!\!\!&&\!\!\!
+\alpha^{M\!L}_{-\frac32,-1}(M_{-\frac52}\otimes L_{1}
-2M_{-\frac32}\otimes L_{0}+M_{-\frac12}\otimes L_{-1}).\label{201303310901}
\end{eqnarray}

The identity $[G_{\frac32},G_{-\frac12}]=2L_1$, together with \eqref{201303300803} and \eqref{201303300810}, gives $G_{-\frac12}\ast d_0(G_{\frac32})=0$. Comparing the relative coefficients, we can deduce the following identities:
\begin{eqnarray}
\!\!\!\!\!\!\!\!
\alpha^{L\!Y}_{\frac32,i}\!\!\!&=&\!\!\!\alpha^{G\!Y}_{\frac32,s}
=\alpha^{Y\!L}_{\frac32,i}=\alpha^{Y\!G}_{\frac32,s}
=\alpha^{Y\!Y}_{\frac32,p}
=\alpha^{Y\!M}_{\frac32,p}
=\alpha^{M\!Y}_{\frac32,p}\nonumber\\
\!\!\!&=&\!\!\!\alpha^{L\!M}_{\frac32,i_1}=\alpha^{G\!M}_{\frac32,s_1}
=\alpha^{M\!L}_{\frac32,i_1}=\alpha^{M\!G}_{\frac32,s_2}
=\alpha^{M\!M}_{\frac32,p_1}=0,\label{201303310800}
\end{eqnarray}
for all $i\in\mathbb{Z}$, $i_1\in\mathbb{Z}\backslash\{\pm1,\,0\}$, $s\in\frac12+\mathbb{Z}$, $s_1\in\frac12+\mathbb{Z}\backslash\{\pm\frac12,\,\frac32\}$, $s_2\in\frac12+\mathbb{Z}\backslash\{-\frac32,\,\pm\frac12\}$, $p\in\frac12\mathbb{Z}$, $p_1\in\frac12\mathbb{Z}\backslash\{-1,\,0,\,-\frac32,\,-\frac12\}$. During the process of comparing the relative coefficients of $G_{-\frac12}\ast d_0(G_{\frac32})=0$, we also can obtain the following identities:
\begin{eqnarray}
\begin{array}{llllllll}\label{201303310900}
&&\alpha^{M\!M}_{\frac32,-\frac32}=-\alpha^{M\!M}_{\frac32,-1}
=-\alpha^{M\!M}_{\frac32,-\frac12}=\alpha^{M\!M}_{\frac32,0},\\[8pt]
&&2\alpha^{L\!M}_{\frac32,-1}=-\alpha^{L\!M}_{\frac32,0}
=2\alpha^{L\!M}_{\frac32,1},\ \ \ 2\alpha^{G\!M}_{\frac32,-\frac12}
=-\alpha^{G\!M}_{\frac32,\frac12}=2\alpha^{G\!M}_{\frac32,\frac32},\\[8pt]
&&2\alpha^{M\!L}_{\frac32,-1}=-\alpha^{M\!L}_{\frac32,0}=
2\alpha^{M\!L}_{\frac32,1},\ \ \ 2\alpha^{M\!G}_{\frac32,-\frac32}=-\alpha^{M\!G}_{\frac32,-\frac12}
=2\alpha^{M\!G}_{\frac32,\frac12}.
\end{array}
\end{eqnarray}

Combining \eqref{201303310600}, \eqref{201303310601} and Convention \ref{201303290600}, we can rewrite $d_0(G_{\frac32})$ as follows:
\begin{eqnarray}
d_0(G_{\frac32})\!\!\!&=&\!\!\!
\alpha^{L\!M}_{\frac32,-1}(L_{-1}\otimes M_{\frac52}
-2L_{0}\otimes M_{\frac32}+L_{1}\otimes M_{\frac12})\nonumber\\
\!\!\!&&\!\!\!
+\alpha^{M\!L}_{\frac32,-1}(M_{\frac12}\otimes L_{1}
-2M_{\frac32}\otimes L_{0}+M_{\frac52}\otimes L_{-1}).\label{201303310906}
\end{eqnarray}

The identity $[G_{\frac32},G_{-\frac32}]=2L_0$, together with Claim \ref{201302182222}, gives $G_{\frac32}\ast d_0(G_{-\frac32})\equiv-G_{-\frac32}\ast d_0(G_{\frac32})$. Using  \eqref{201303310901} and \eqref{201303310906} and comparing the relative coefficients, we can obtain the following identities:
\begin{eqnarray*}
(i+6)\alpha^{M\!L}_{-\frac32,i+3}+(i-3)\alpha^{M\!L}_{\frac32,i}\!\!\!&=&\!\!\!
(i-3)\alpha^{L\!M}_{-\frac32,i}+(i+6)\alpha^{L\!M}_{\frac32,i+3}\\
\!\!\!&=&\!\!\!\alpha^{L\!M}_{-\frac32,i}+\alpha^{L\!M}_{\frac32,i}=
\alpha^{M\!L}_{-\frac32,i}+\alpha^{M\!L}_{\frac32,i}=0,\ \ \ \forall\,\,i\in\mathbb{Z},
\end{eqnarray*}
which together with \eqref{201303310901} and \eqref{201303310906}, force
\begin{eqnarray*}
\alpha^{L\!M}_{-\frac32,i}\!\!\!&=&\!\!\!\alpha^{L\!M}_{\frac32,i}
=\alpha^{M\!L}_{-\frac32,i}=\alpha^{M\!L}_{\frac32,i}=0,\ \ \ \forall\,\,i\in\mathbb{Z}.
\end{eqnarray*}
Then $d_0(-G_{\frac32})$ and $d_0(G_{\frac32})$ referred in \eqref{201303310901} and \eqref{201303310906} can be simplified as follows:
\begin{eqnarray}
d_0(-G_{\frac32})\!\!\!&=&\!\!\!d_0(G_{\frac32})=0.\label{201303310909}
\end{eqnarray}

For the given $d_0\in{\rm Der}(\frak{tsns},\mathcal{H}/\frak{tsns}_{\mathcal{C}})$, we can write $d_0(Y_{\pm\frac12})$ as follows:
\begin{eqnarray*}
&&\!\!\!\!\!\!\!\!d_0(Y_{\pm\frac12})=
\mbox{$\sum\limits_{i\in\mathbb{Z}}$}
(\beta^{L\!Y}_{\pm\frac12,i}L_{i}\otimes Y_{\pm\frac12-i}
+\beta^{Y\!L}_{\pm\frac12,i}Y_{i\pm\frac12}\otimes L_{-i}
+\beta^{L\!M}_{\pm\frac12,i}L_{i}\otimes M_{\pm\frac12-i}
+\beta^{M\!L}_{\pm\frac12,i}M_{i\pm\frac12}\otimes L_{-i})\\
&&\!\!\!\!\!\!\!\!
+\mbox{$\sum\limits_{s\in\frac{1}{2}+\mathbb{Z}}$}
(\beta^{G\!Y}_{\pm\frac12,s}G_{s}\otimes Y_{\pm\frac12-s}
+\beta^{Y\!G}_{\pm\frac12,s}Y_{s\pm\frac12}\otimes G_{-s}
+\beta^{G\!M}_{\pm\frac12,s}G_{s}\otimes M_{\pm\frac12-s}
+\beta^{M\!G}_{\pm\frac12,s}M_{s\pm\frac12}\otimes G_{-s})\\
&&\!\!\!\!\!\!\!\!
+\mbox{$\sum\limits_{p\in\frac12\mathbb{Z}}$}
(\beta^{Y\!Y}_{\pm\frac12,p}Y_{p\pm\frac12}\otimes Y_{-p}
+\beta^{Y\!M}_{\pm\frac12,p}Y_{p\pm\frac12}\otimes M_{-p}
+\beta^{M\!Y}_{\pm\frac12,p}M_{p\pm\frac12}\otimes Y_{-p}
+\beta^{M\!M}_{\pm\frac12,p}M_{p\pm\frac12}\otimes M_{-p}),
\end{eqnarray*}
where the coefficients are all in $\mathbb{C}$ and the sums are all finite. The identities $[Y_{\frac12},M_{-\frac12}]=0$ and $[Y_{-\frac12},M_{\frac12}]=0$, together with Claim \ref{201302182222}, \eqref{201303300816} and \eqref{201303300906}, gives $M_{-\frac12}\ast d_0(Y_{\frac12})\equiv M_{\frac12}\ast d_0(Y_{-\frac12})\equiv0$. Comparing the relative coefficients, we can obtain the following identities:
\begin{eqnarray}
\beta^{G\!Y}_{\frac12,s}\!\!\!&=&\!\!\!\beta^{Y\!G}_{\frac12,s}=
(1-i)\beta^{L\!Y}_{\frac12,i}=(i+1)\beta^{Y\!L}_{\frac12,i}\nonumber\\
\!\!\!&=&\!\!\!\beta^{G\!M}_{\frac12,s+1}+\beta^{M\!G}_{\frac12,s}=
i\beta^{L\!M}_{\frac12,i+1}+(i+1)\beta^{M\!L}_{\frac12,i}=0,\label{201304010600}\\
\beta^{G\!Y}_{-\frac12,s}\!\!\!&=&\!\!\!\beta^{Y\!G}_{-\frac12,s}
=(i+1)\beta^{L\!Y}_{-\frac12,i}=(1-i)\beta^{Y\!L}_{-\frac12,i}\nonumber\\
\!\!\!&=&\!\!\!\beta^{G\!M}_{-\frac12,s}+\beta^{M\!G}_{-\frac12,s+1}=
i\beta^{M\!L}_{-\frac12,i+1}+(i+1)\beta^{L\!M}_{-\frac12,i}=0,\label{201304010601}
\end{eqnarray}
for all $i\in\mathbb{Z}$, $s\in\frac12+\mathbb{Z}$.

The identities $[G_{\frac32},M_{-\frac12}]=2M_1=[Y_{\frac12},Y_{\frac12}]$ and $[G_{-\frac32},M_{\frac12}]=2M_{-1}=[Y_{-\frac12},Y_{-\frac12}]$, together with the obtained results $d_0(G_{\pm\frac32})=d_0(M_{\pm\frac12})=0$, give $Y_{\frac12}\ast d_0(Y_{\frac12})=Y_{-\frac12}\ast d_0(Y_{-\frac12})=0$. Comparing the relative coefficients, we can obtain the following identities:
\begin{eqnarray}
\beta^{L\!Y}_{\frac12,i}\!\!\!&=&\!\!\!\beta^{Y\!L}_{\frac12,i}
=2\beta^{L\!M}_{\frac12,i}-\beta^{Y\!Y}_{\frac12,i}
=2\beta^{M\!L}_{\frac12,i}-\beta^{Y\!Y}_{\frac12,i-\frac12}
=\beta^{G\!M}_{\frac12,s}+\beta^{Y\!Y}_{\frac12,s}\nonumber\\
\!\!\!&=&\!\!\!\beta^{M\!G}_{\frac12,s}+\beta^{Y\!Y}_{\frac12,s-\frac12}
=\beta^{Y\!M}_{\frac12,i}+\beta^{M\!Y}_{\frac12,i+\frac12}
=\beta^{Y\!M}_{\frac12,i-\frac12}-\beta^{M\!Y}_{\frac12,i}=0,\label{201304010602}\\
\beta^{L\!Y}_{-\frac12,i}\!\!\!&=&\!\!\!\beta^{Y\!L}_{-\frac12,i}
=2\beta^{L\!M}_{-\frac12,i}-\beta^{Y\!Y}_{-\frac12,i}
=2\beta^{M\!L}_{-\frac12,i}-\beta^{Y\!Y}_{-\frac12,i+\frac12}
=\beta^{G\!M}_{-\frac12,s}+\beta^{Y\!Y}_{-\frac12,s}\nonumber\\
\!\!\!&=&\!\!\!\beta^{M\!G}_{-\frac12,s}+\beta^{Y\!Y}_{-\frac12,s+\frac12}
=\beta^{Y\!M}_{-\frac12,i}+\beta^{M\!Y}_{-\frac12,i-\frac12}
=\beta^{M\!Y}_{-\frac12,i}-\beta^{Y\!M}_{-\frac12,i+\frac12}=0,\label{201304010603}
\end{eqnarray}
for all $i\in\mathbb{Z}$,  $s\in\frac12+\mathbb{Z}$. Combining the identities gotten in \eqref{201304010600}, \eqref{201304010601}, \eqref{201304010602} and \eqref{201304010603}, we can deduce the following identities:
\begin{eqnarray}
\beta^{L\!Y}_{\frac12,i}\!\!\!&=&\!\!\!\beta^{Y\!L}_{\frac12,i}
=\beta^{L\!M}_{\frac12,i}=\beta^{M\!L}_{\frac12,i}
=\beta^{G\!Y}_{\frac12,s}=\beta^{Y\!G}_{\frac12,s}
=\beta^{G\!M}_{\frac12,s}=\beta^{M\!G}_{\frac12,s}
=\beta^{Y\!Y}_{\frac12,p}=0,\label{201304010605}\\
\beta^{L\!Y}_{-\frac12,i}\!\!\!&=&\!\!\!\beta^{Y\!L}_{-\frac12,i}
=\beta^{L\!M}_{-\frac12,i}=\beta^{M\!L}_{-\frac12,i}
=\beta^{G\!Y}_{-\frac12,s}=\beta^{Y\!G}_{-\frac12,s}
=\beta^{G\!M}_{-\frac12,s}=\beta^{M\!G}_{-\frac12,s}
=\beta^{Y\!Y}_{-\frac12,p}=0,\label{201304010606}
\end{eqnarray}
for all $i\in\mathbb{Z}$,  $s\in\frac12+\mathbb{Z}$ and  $p\in\frac12\mathbb{Z}$. Then $d_0(Y_{\pm\frac12})$ can be rewritten as follows:
\begin{eqnarray}
d_0(Y_{\pm\frac12})\!\!\!&=&\!\!\!
\mbox{$\sum\limits_{p\in\frac12\mathbb{Z}}$}
(\beta^{Y\!Y}_{\pm\frac12,p}Y_{p\pm\frac12}\otimes Y_{-p}
+\beta^{Y\!M}_{\pm\frac12,p}Y_{p\pm\frac12}\otimes M_{-p}\nonumber\\
\!\!\!&&\!\!\!
+\beta^{M\!Y}_{\pm\frac12,p}M_{p\pm\frac12}\otimes Y_{-p}
+\beta^{M\!M}_{\pm\frac12,p}M_{p\pm\frac12}\otimes M_{-p}).\label{201304010608}
\end{eqnarray}
Up to now, we have obtained the following results:
\begin{eqnarray}\label{201304010609}
d_0(L_{1})\!\!\!&=&\!\!\!d_0(G_{\pm\frac12})
=d_0(G_{\pm\frac32})=d_0(M_{\pm\frac12})=0.
\end{eqnarray}

The following identities have been obtained in \eqref{201304010602} and \eqref{201304010603}:
\begin{eqnarray}
&&\beta^{Y\!M}_{\frac12,i}+\beta^{M\!Y}_{\frac12,i+\frac12}
=\beta^{Y\!M}_{\frac12,i-\frac12}-\beta^{M\!Y}_{\frac12,i}
=\beta^{Y\!M}_{-\frac12,i}+\beta^{M\!Y}_{-\frac12,i-\frac12}
=\beta^{M\!Y}_{-\frac12,i}-\beta^{Y\!M}_{-\frac12,i+\frac12}=0,\label{201304010610}
\end{eqnarray}
for all $i\in\mathbb{Z}$,  $s\in\frac12+\mathbb{Z}$.

The identities $[L_1,Y_{-\frac12}]=-\frac12Y_{\frac12}$, together with \eqref{201304010609}, gives
\begin{eqnarray}\label{201304011700}
&&2L_1\ast d_0(Y_{-\frac12})+d_0(Y_{\frac12})=0.
\end{eqnarray}
Comparing the relative coefficients, we can obtain the following identities:
\begin{eqnarray}
\begin{array}{lll}\label{201304011200}
2i\beta^{M\!Y}_{-\frac12,i+\frac12}-(2i+3)\beta^{M\!Y}_{-\frac12,i+\frac32}
+\beta^{M\!Y}_{\frac12,i+\frac12}=0,\\[8pt]
(2i-1)\beta^{Y\!M}_{-\frac12,i}-2(i+1)\beta^{Y\!M}_{-\frac12,i+1}
+\beta^{Y\!M}_{\frac12,i}=0,
\end{array}
\end{eqnarray}
for all $i\in\mathbb{Z}$, from which we can deduce
\begin{eqnarray}\label{201304010611}
&&2(2i+1)\beta^{Y\!M}_{-\frac12,i+1}
=(2i+3)\beta^{Y\!M}_{-\frac12,i+2}+(2i-1)\beta^{Y\!M}_{-\frac12,i},\ \ \ \forall\,\,i\in\mathbb{Z}.
\end{eqnarray}
According to Convention \ref{201303290600}, we know that $\beta^{Y\!M}_{-\frac12,0}=0$. Then $\beta^{Y\!M}_{-\frac12,i}$ can be linearly expressed by $\beta^{Y\!M}_{-\frac12,1}$ for any $i\geqslant1$ while  $\beta^{Y\!M}_{-\frac12,j}$ can be linearly expressed by $\beta^{Y\!M}_{-\frac12,-1}$ for any $j\leqslant-1$. Then the following identity follows from the finiteness of the relative sum referred in \eqref{201304010608}:
\begin{eqnarray}\label{201304011000}
\beta^{Y\!M}_{-\frac12,i}=0,\ \ \ \forall\,\,i\in\mathbb{Z},
\end{eqnarray}
And the following identities follows from \eqref{201304010610}, \eqref{201304011200} and \eqref{201304011000}:
\begin{eqnarray}\label{201304011001}
&&\beta^{M\!Y}_{-\frac12,s}=\beta^{Y\!M}_{\frac12,i}
=\beta^{M\!Y}_{\frac12,s}=0
\end{eqnarray}

Comparing the relative coefficients of $2L_1\ast d_0(Y_{-\frac12})+d_0(Y_{\frac12})=0$ referred in \eqref{201304011700}, we can obtain the following identities:
\begin{eqnarray}\label{201304011800}
&&2(s-1)\beta^{Y\!M}_{-\frac12,s}-(2s+1)\beta^{Y\!M}_{-\frac12,s+1}
+\beta^{Y\!M}_{\frac12,s}=0,\ \ \ \forall\,\,s\in\frac12+\mathbb{Z}.
\end{eqnarray}
The identities $[G_{\frac12},Y_{\frac12}]=2Y_1
=[G_{\frac32},Y_{-\frac12}]$, together with \eqref{201304010609}, gives
$G_{\frac12}\ast d_0(Y_{\frac12})=G_{\frac32}\ast d_0(Y_{-\frac12})$.
Comparing the coefficients of $Y_{i+1}\otimes M_{-i}$, we can obtain the following identities:
\begin{eqnarray*}
&&\beta^{Y\!M}_{\frac12,i}+\beta^{Y\!M}_{\frac12,i+\frac12}
=\beta^{Y\!M}_{-\frac12,i}+\beta^{Y\!M}_{-\frac12,i+\frac32},
\end{eqnarray*}
for all $i\in\mathbb{Z}$, which together with \eqref{201304011000} and \eqref{201304011001}, forces
\begin{eqnarray}\label{201304011801}
&&\beta^{Y\!M}_{\frac12,s}=\beta^{Y\!M}_{-\frac12,s+1},
\ \ \ \forall\,\,s\in\frac12+\mathbb{Z}.
\end{eqnarray}
Then combining \eqref{201304011800} and \eqref{201304011801}, we can obtain
\begin{eqnarray*}
&&(s-1)\beta^{Y\!M}_{-\frac12,s}=s\beta^{Y\!M}_{-\frac12,s+1},
\ \ \ \forall\,\,s\in\frac12+\mathbb{Z},
\end{eqnarray*}
which forces
\begin{eqnarray}\label{201304011802}
&&\beta^{Y\!M}_{-\frac12,s}=0,\ \ \ \forall\,\,s\in\frac12+\mathbb{Z}.
\end{eqnarray}
According to \eqref{201304010610}, \eqref{201304011800} and \eqref{201304011802}, we can deduce the following identities:
\begin{eqnarray}\label{201304011803}
&&\beta^{Y\!M}_{-\frac12,s}=\beta^{Y\!M}_{\frac12,s}
=\beta^{M\!Y}_{\frac12,i}=\beta^{M\!Y}_{-\frac12,i}=0,
\ \ \ \forall\,\,i\in\mathbb{Z},\,s\in\frac12+\mathbb{Z}.
\end{eqnarray}
According to \eqref{201304011000}, \eqref{201304011001}, \eqref{201304011802}, \eqref{201304011803}, we can rewrite \eqref{201304010608} as follows:
\begin{eqnarray}\label{201304012000}
d_0(Y_{\pm\frac12})\!\!\!&=&\!\!\!
\mbox{$\sum\limits_{p\in\frac12\mathbb{Z}}$}
\beta^{M\!M}_{\pm\frac12,p}M_{p\pm\frac12}\otimes M_{-p}.
\end{eqnarray}

The identities $[G_{-\frac12},[G_{\frac12},Y_{\frac12}]]=2[G_{-\frac12},Y_1]=\frac32Y_{\frac12}$, together with \eqref{201304010609}, gives
$G_{-\frac12}\ast G_{\frac12}\ast d_0(Y_{\frac12})=\frac32d_0(Y_{\frac12})$.
Comparing the relative coefficients, we can obtain the following:
\begin{eqnarray}
&&(i+1)\beta^{M\!M}_{\frac12,i+1}-i\beta^{M\!M}_{\frac12,i}
=\frac12\beta^{M\!M}_{\frac12,i+\frac12},\label{201304012001}\\
&&(i+1)\beta^{M\!M}_{\frac12,i+\frac12}-i\beta^{M\!M}_{\frac12,i-\frac12}
=\frac12\beta^{M\!M}_{\frac12,i},\label{201304012002}
\end{eqnarray}
for all $i\in\mathbb{Z}$. Then we can deduce
\begin{eqnarray*}
&&2(i+1)\big((i+1)\beta^{M\!M}_{\frac12,i+1}-i\beta^{M\!M}_{\frac12,i}\big)
-2i\big(i\beta^{M\!M}_{\frac12,i}-(i-1)\beta^{M\!M}_{\frac12,i-1}\big)
=\frac12\beta^{M\!M}_{\frac12,i},
\end{eqnarray*}
for all $i\in\mathbb{Z}$, which forces (noticing $\beta^{M\!M}_{\frac12,0}=0$ according to Convention \ref{201303290600})
\begin{eqnarray}\label{201304012006}
&&\beta^{M\!M}_{\frac12,i}=0,\ \ \ \forall\,\,i\in\mathbb{Z}.
\end{eqnarray}
Recalling  \eqref{201304012001}, we can deduce
\begin{eqnarray}\label{201304012008}
&&\beta^{M\!M}_{\frac12,s}=0,\ \ \ \forall\,\,s\in\frac12+\mathbb{Z}.
\end{eqnarray}

Comparing the relative coefficients of $2L_1\ast d_0(Y_{-\frac12})+d_0(Y_{\frac12})=0$ referred in \eqref{201304011700}, we can obtain the following identities:
\begin{eqnarray}\label{201304012009}
\begin{array}{lll}
&&2i\beta^{M\!M}_{-\frac12,i}-2(i+1)\beta^{M\!M}_{-\frac12,i+1}
+\beta^{M\!M}_{\frac12,i}=0,\\[8pt]
&&(2s-1)\beta^{M\!M}_{-\frac12,s}-(2s+1)\beta^{M\!M}_{-\frac12,s+1}
+\beta^{M\!M}_{\frac12,s}=0,
\end{array}
\end{eqnarray}
for all $i\in\mathbb{Z}$, $s\in\frac12+\mathbb{Z}$. Combining \eqref{201304012006}, \eqref{201304012008}, \eqref{201304012009} and recalling $\beta^{M\!M}_{-\frac12,0}=\beta^{M\!M}_{-\frac12,\frac12}=0$ according to Convention \ref{201303290600}, we can deduce
\begin{eqnarray}\label{201304012006}
&&\beta^{M\!M}_{-\frac12,p}=0,\ \ \ \forall\,\,p\in\frac12\mathbb{Z}.
\end{eqnarray}
Hence \eqref{201304012000} can be simplified as follows:
\begin{eqnarray}\label{201304012200}
d_0(Y_{\pm\frac12})\!\!\!&=&\!\!\!0.
\end{eqnarray}

Then Claim \ref{201303091203} finally follows from \eqref{201304010609} and \eqref{201304012200}.

Thus $H^1(\frak{tsns},\frak{tsns}_{\mathcal{C}})\cong H^1(\frak{tsns},\mathcal{H})$ according to Claim \ref{201303091203} and \eqref{201303100606}.

\begin{clai}\label{201303091209}
$H^1(\frak{tsns},\frak{tsns}_{\mathcal{C}})\cong\frak{D}^\natural$.
\end{clai}
This claim follows from Lemmas \ref{201303070600} and \ref{201303091802}.

\begin{clai}\label{201303091201}
For any $d\in{\rm Der}(\frak{tsns},\frak{tsns}^\otimes)$, \eqref{201303091200} is a finite sum.
\end{clai}

For any $p\in\frac12\mathbb{Z}$, one can suppose $d_p=(v_p)_{\rm inn}+\delta_{p,0}\mathfrak{d}^\natural$ for some $v_p\in\frak{tsns}^\otimes_p$ and $\mathfrak{d}^\natural\in\frak{D}^\natural$. If
$\Delta=\{p\in\frac12\mathbb{Z}\,|\,v_p\neq0\}$ is an infinite set, then
$d(L_0)=\mathfrak{d}^\natural(L_0)+\sum_{p\in\Delta}L_0\ast
v_p=\mathfrak{d}^\natural(L_0)-\sum_{p\in\Delta}p\,v_p$ is an infinite sum, which is not in $\frak{tsns}^\otimes$, contradicting the fact that $d$ is a derivation from $\frak{tsns}$ to $\frak{tsns}^\otimes$. Then Claim \ref{201303091201} follows.

By now we have completed the proof of Theorem \ref{201303090909}\,(i).\hfill$\Box$\par\vskip6pt

The following lemma is still true for $\frak{tsns}$ by employing the technique of \cite[Lemma 3.5]{YS--AMS2010} or  \cite[Lemma 2.5]{FLX--AC2011}.

\begin{lemm}\label{201303091600}
Suppose $v\in\frak{tsns}^\otimes\big({\rm modulo}\,(\mathcal{C}^\otimes)\big)$ such that $x\cdot v\in{\rm Im}({\bf\it1}\otimes{\bf\it1}-\tau)$ for all
$x\in\frak{tsns}.$ Then there exists some $u\in{\rm Im}({\bf\it1}\otimes{\bf\it1}-\tau)$ such that $v-u\in\mathcal{C}^\otimes$.
\end{lemm}
\noindent{\it Proof of Theorem \ref{201303090909}\,(ii),\,(iii)}\ \ \,They follow from Lemmas \ref{201303070606}, \ref{201303091600} and Theorem \ref{201303090909}\,(i).\hfill$\Box$\par

\vspace*{12pt}

\noindent{\bf Acknowledgements}\ \,This work was supported by a NSF grant BK20160403 of Jiangsu Province and NSF grants 11671056, 11271056, 11101056 of China.

\end{document}